\documentclass[11pt]{article}
\usepackage{float}
\usepackage{amsfonts,amsmath,amssymb,graphicx,color}
\usepackage{amsthm}
\usepackage{algorithm}
\usepackage{algpseudocode}
\usepackage{psfrag}
\usepackage{pdflscape}
\usepackage{multirow} 
\usepackage{enumitem}
\usepackage{framed}
\usepackage{sidecap}
\usepackage{tikz}
\usepackage{subcaption}
\usepackage{lscape}
\usepackage[scientific-notation=true]{siunitx}
\usetikzlibrary{arrows}
\usetikzlibrary{patterns,arrows,decorations.pathreplacing}
\graphicspath{ {./figures/} }
    \usepackage{soul}
\usepackage[normalem]{ulem}
\usepackage{url}
\usepackage[numbers, sort&compress ]{natbib} 
\usepackage{microtype}
\usepackage{hyperref}

\numberwithin{figure}{section}
\numberwithin{table}{section}

\newcommand{\R}{\mathbb{R}}

\numberwithin{equation}{section}

\newcommand{\comment}[1]{}

\newcommand{\be}{\begin{equation}}
\newcommand{\ee}{\end{equation}}

\newcommand{\bea}{\begin{eqnarray}}
\newcommand{\eea}{\end{eqnarray}}
\newcommand{\beqa}{\begin{eqnarray}}
\newcommand{\eeqa}{\end{eqnarray}}
\newcommand{\beann}{\begin{eqnarray*}}
\newcommand{\eeann}{\end{eqnarray*}}

\usepackage{scalerel}

\newcommand{\bmat}{\left[ \begin{array}}
\newcommand{\emat}{\end{array} \right]}

\newcommand{\beq}{\begin{equation}}
\newcommand{\eeq}{\end{equation}}

\newcommand{\bproof}{\begin{description} \item[{\it Proof}.] ~ }
\newcommand{\eproof}{\hspace*{\fill}$\Box$\medskip \end{description}}

\newcounter{algo}[section]

\newcounter{prog}[section]

\title{A Feasible  Method for Constrained  Derivative-Free Optimization}
\author{       
        Melody Qiming Xuan\thanks{Department of Industrial Engineering and Management Sciences, Northwestern University,  USA. This author was supported by AFOSR grant FA95502110084 and ONR grant N00014-21-1-2675}
               \and 
       Jorge Nocedal \thanks{Department of Industrial Engineering and Management Sciences, Northwestern University,  USA.  This author was supported by National Science Foundation grant DMS-2011494, AFOSR grant FA95502110084, and  ONR grant N00014-21-1-2675.} 
      }
\date{February 16, 2024}
\begin{document}

\maketitle

\begin{abstract}
This paper explores a method for solving constrained optimization problems when the derivatives of the objective function are unavailable, while the derivatives of the constraints are known. We allow the objective and constraint function to be nonconvex.
The method constructs a quadratic model of the objective function via interpolation and  computes a step by minimizing this model subject to the original constraints in the problem and a trust region constraint. The step computation requires the solution of a general nonlinear program, which is economically feasible when the constraints and their derivatives are very inexpensive to compute compared to the objective function. The paper includes a summary of numerical results that highlight the method's promising potential.

\end{abstract}

%\tableofcontents

\bigskip
%\jn{Jorge}, \mx{Melody}
\section{Introduction}
\label{intro}
\setcounter{equation}{0}

Interpolation-based trust-region optimization \texttt{(IBO)} methods are among the most efficient and reliable techniques for unconstrained derivative-free optimization (DFO); see e.g. \cite{more2009benchmarking,rios2013derivative,shi2022numerical}. Yet extending \texttt{IBO} methods to handle general constraints is not straightforward and has received limited attention in the literature \cite{larson2019derivative}. In this paper, we investigate a method first proposed by Conn et al. \cite{conn1998derivative} that is designed for problems where the evaluation of the constraints is inexpensive compared to the evaluation of the objective function---a common setting in DFO applications. The method described in \cite{conn1998derivative} is a feasible method that constructs a quadratic model of the objective using an IBO approach and computes a step by minimizing this model subject to the original constraints of the problem and a trust region constraint.  Thus, in general, the step computation requires the solution of a nonlinear optimization subproblem, which can be done using a general purpose nonlinear programming method. We refer to this approach as \texttt{FIBO}, for feasible interpolation-based optimization method. 

The numerical results reported in \cite{conn1998derivative}, and subsequently in \cite{conejo2015trust}, are, however, inconclusive
and the \texttt{FIBO} approach is rarely (if ever) used in practice. In this paper, we give careful consideration to the design of \texttt{FIBO} and present numerical results that indicate that \texttt{FIBO} is one of the most efficient methods for important classes of constrained DFO applications. An added appeal of \texttt{FIBO} is that it can be easily constructed around an unconstrained \texttt{IBO} method, and can use many nonlinear programming solvers to compute the algorithm's step. This sets it apart from other methods proposed for constrained DFO applications,
which rely on heuristic penalty or augmented Lagrangian approaches \cite{lewis2002globally,griffin2010nonlinearly}
or involve complex \texttt{IBO} methodology \cite{lincoa}. Our work goes beyond \cite{hough2022model} in that we allow the constraints to be nonconvex.

The problem under consideration is
\begin{align}  \min_x & \ \ f(x) \nonumber \\
  \mbox{s.t.} & \  \ c_i(x) = 0, \quad i \in {\mathcal E}, \label{problem}
\\
                    & \  \ c_i(x) \leq 0, \quad i \in {\mathcal I}, \nonumber
\end{align}
for some finite index sets $\mathcal E$ and $\mathcal I$, and
where $f: \R^n \rightarrow \R$ and  $c_i(x) : \R^n \rightarrow \R, \, \, i  \in \mathcal{I} \cup \mathcal{E},$ are smooth functions.  We assume that the derivatives of $f$ are not available. We also assume---and this is important---that the evaluation of the functions $c_i$ and their derivatives is quite inexpensive compared to the cost of evaluating $f$.

\section{The FIBO Algorithm}
Let us begin by reviewing the basic \texttt{IBO} method for \emph{unconstrained} derivative-free optimization  \cite{conn2009introduction,powell2006newuoa}, which has proved to be very effective in practice. The algorithm starts by evaluating the objective function at a set of interpolation points $Y_0$, normally placed on a stencil around the initial point $x_0$. The algorithm then sets $ k \leftarrow 0$ and, at every iteration, constructs a quadratic model of the objective function in the form
\begin{equation}\label{eq:quad}
	m_k(x_k + s) = f(x_k) + g_k^Ts + \frac{1}{2}s^TH_ks.
\end{equation}
To define this model, we need to specify the vector $g_k$ and the symmetric matrix $H_k$. 
This  can be defined in various ways. Assuming that $Y_k$ is \textit{poised} \cite{conn2009introduction} and that $|Y_k| = (n+1)(n+2)/2$, the model $m_k$ is uniquely specified by requiring it to interpolate the objective function $f$  at the current interpolation set $Y_k$. This is the vanilla \texttt{IBO} method, and requires $O(n^6)$ work per iteration. A more sophisticated \texttt{IBO} approach proposed by Powell \cite{powell2006newuoa} defines $m_k$ through a minimum Frobenius norm update of the matrix $H_k$ subject to $|Y_k|= 2n +1$ interpolation conditions. By updating matrix factorizations and reusing information judiciously, the computational cost of the iteration is dramatically reduced to $O(n^2)$ operations. The method proposed in  \cite{powell2006newuoa} provides great flexibility in that the size of interpolation set $Y_k$ can be chosen from $n+2$ to $(n+1)(n+2)/2$ points (with the choice $2n+1$ being the default).

Regardless of how the model is constructed, the trial step $s_k$ of \texttt{IBO} methods for unconstrained optimization is given by the solution of the trust region subproblem 
\begin{subequations}
\label{eq:tr sub}
\begin{align} 
\min_s  &  \quad m_k(x_k+s)  \label{qmodel} \\ 
 \mbox{s.t.} & \quad \| s\|_2 \leq \Delta_k. 
\end{align}
\end{subequations}
As in any trust region method, acceptance of the step $s_k$ is  determined based on the ratio
\begin{align}\label{eq:ratio test}
    \rho_k = \frac{f(x_k) - f(x_k + s_k)}{m_k(x_k) - m_k(x_k + s_k)}
\end{align}
between actual vs predicted reduction in the objective. If sufficient decrease is obtained (i.e., $\rho_k \geq \eta$ for some $\eta \in (0,1)$), the step is considered successful and the new iterate is defined as $x_{k+1} = x_k + s_k$; otherwise the step is rejected and  $x_{k+1} \gets x_k$. 

The trust region radius $\Delta_k$ is increased whenever the step is successful. Otherwise, the algorithm considers two cases. If the interpolation set $Y_k$ is not poised, (i.e., if the interpolation points lie (nearly) in a linear subspace of $\mathbb{R}^n$), then $\Delta_k$ remains unchanged to avoid premature shrinkage of the trust region, and $Y_k$ is improved via a \emph{geometry phase} \cite{powell2006newuoa,conn2009introduction,fasano2009geometry}. If, on the other hand, $Y_k$ is poised, then $\Delta_k$ is decreased by a constant factor as in a standard trust region method for derivative-based optimization. 

The algorithm's final component involves the strategy for updating the interpolation set $Y_k$ at each iteration. The most commonly employed approach is to remove the point that is farthest from $x_{k+1}$ and replace it with $x_{k+1}$. Other strategies for updating $Y_k$  have been employed \cite{powell2006newuoa}, but the specific details are not crucial to the ensuing discussion.

\paragraph{Specification of the Algorithm}

We  now describe \texttt{FIBO}, the extension of the \texttt{IBO} method for the solution of 
 the constrained optimization problem \eqref{problem}.  In cases where evaluating the constraints and their derivatives is  inexpensive compared to the objective function, it is realistic to compute the algorithm's step by solving the subproblem.
  \begin{subequations}
\label{outer}
\begin{align} 
\min_s  &  \ \  m_k(x_k + s) = f(x_k) + g_k^Ts  + \tfrac{1}{2} s^T H_k s   \label{modela} \\ 
 \mbox{s.t.} & \  \ c_i(x_k+s) = 0, \quad i \in {\mathcal E}, \label{m1}\\
                    & \  \ c_i(x_k+s) \leq 0, \quad i \in {\mathcal I}, \label{m2}\\
                    & \ \ \| s\|_2 \leq \Delta_k. \label{tr radius}
\end{align}
\end{subequations}
The model $m_k$ in \eqref{modela} is constructed, as before, by interpolating function values at a set of interpolation points $Y_k$. Thus, as in the unconstrained case $m_k$ models the objective function $f$. Note that \eqref{outer} includes the original constraints of the problem (not their linearizations) and therefore \eqref{outer} is in general a (nonconvex) nonlinear programming problem that must be solved using a general-purpose nonlinear programming method (\texttt{KNITRO} in our experiments). All other details of the algorithm, such as the update of the interpolation set and trust region, are defined as in the unconstrained  \texttt{IBO} method. In fact, to develop the \texttt{FIBO} method  one can start with any IBO method for unconstrained optimization and  replace the step computation \eqref{eq:tr sub} by \eqref{outer}.
 We state the proposed method in Algorithm \ref{alg:outerDFO}.

\begin{algorithm}[H]
    \caption{FIBO}
    \begin{algorithmic}[1]
        \State{Parameters: $0 \leq \eta < 1$, and $\Delta_0 >0$. }
        \State{Choose $x_0$, a \textbf{feasible} starting point, and construct an initial  interpolation set $Y_0$.}
        \State{Set $k = 0$.}
        \While{no convergence test is satisfied}  
        \State{Build a local (quadratic) model \eqref{modela} using the interpolation set $Y_k$.}           
        
        \State{Compute a step $s_k$ by solving the trust region subproblem \eqref{outer}.}
        \State{Compute the ratio defined by \eqref{eq:ratio test}.}
        \If{$\rho_k \geq \eta$}
        \State{Set $x_{k+1} = x_{k}+s_k$. }
        \State{Choose $\Delta_{k+1} \geq \Delta_k$}.
        \State{Update $Y_{k}$ to include $x_{k+1}$.}
        \ElsIf{$Y_k$ needs to be improved}
        \State{Set $x_{k+1} = x_k$.}
        \State{Set $\Delta_{k+1} = \Delta_k$. }
        \State{Improve $Y_k$ using a geometry-improving procedure.}
        \Else
        \State{Set $x_{k+1} = x_k$.}
        \State{Choose $\Delta_{k+1} < \Delta_k$.}
        \State{Update $Y_{k}$ to include $x_{k}+s_k$.}
        \EndIf 
        \State{$Y_{k+1} = Y_k$ }
        \State{$k \leftarrow k + 1$}
        \EndWhile
    \end{algorithmic}\label{alg:outerDFO}
    \end{algorithm}

In order to fully characterize the algorithm, one must specify the methodology for updating the interpolation set (lines 11 and 19)  and the procedure for improving its geometry (line 15). We do not provide specifics here since every practical \texttt{IBO} method has its own way of implementing these procedures, which can then be incorporated into \texttt{FIBO}.  In section~\ref{sec:num}, we describe the methodology employed in our implementation.

We now comment on the main properties of the \texttt{FIBO} approach.

\paragraph{\it Feasibility} Under the assumption that the nonlinear programming solver is successful in every iteration (line 6), all iterates of the \texttt{FIBO} algorithm maintain feasibility. This property is highly desirable for DFO problems wherein the objective function is undefined outside the feasible set. (In this situation, the initial interpolation set $Y_0$ must consist only of feasible points.)
It is worth noting that while the nonlinear programming method employed to solve the subproblem \eqref{outer} may generate infeasible iterates, this method only evaluates the quadratic objective $m_k$ of this subproblem, rather than the objective function $f$. Note also that since all iterates are feasible, it is  appropriate for $m_k$ to model only the objective and disregard the contribution of the constraints.

\paragraph{\it Per-Iteration Cost}
The objective function $f$ is evaluated once per iteration.  The cost of constructing the quadratic model depends on the \texttt{IBO} approach used. The vanilla \texttt{IBO} method computes  the model $m_k$ from scratch at every iteration, at the cost of $O(n^6)$ flops.  In contrast, Powell's Hessian update approach, as implemented in \texttt{NEWUOA}, requires only $O(n^2)$ flops to form the model.
   
  The cost of solving the subproblem \eqref{outer} is a crucial factor in determining the feasibility of the \texttt{FIBO} approach. The nonlinear solver may require a significant number of iterations to compute a solution of \eqref{outer} to reasonable accuracy if the constraints \eqref{m1}, \eqref{m2} are not simple. Thus, \texttt{FIBO}  may perform a large number of evaluation of these constraints, as well as many evaluations of the quadratic model $m_k$, the latter at a cost of $O(n^2)$ flops since the Hessian $H_k$ of $m_k$ is dense.             
%    When the number of variables $n$ is very small (say less than 100), the cost of solving \eqref{outer} is of no concern even if the nonlinear programming algorithm performs many iterations to return an accurate solution, \textit{provided} the constraints are inexpensive to evaluation. But as the number of variables becomes large the total linear algebra cost on the nonlinear programming is of concern.  Similarly, if the cost of evaluating the constraints and their derivatives is high, then the \texttt{FIBO} approach is not viable given the need for multiple evaluations for the solution of a single trust region step. In the next section, we quantify more precisely the class of problems for which \texttt{FIBO} is efficient.

\medskip\noindent{\it Robustness.}  The success of the FIBO approach hinges on the ability of the nonlinear programming method to compute an optimal solution of \eqref{outer}---or at least a feasible solution  that makes sufficient progress.
%, as outlined in the convergence analysis discussed in the next section.
While modern nonlinear solvers are quite reliable, there is  no guarantee of their  success in solving problem \eqref{outer}. Nevertheless, it is important to recognize that attaining robustness poses a challenge for all methods designed for constrained DFO.

\medskip\noindent{\it Application to Noisy Problems} The \texttt{FIBO} method can be applied to noisy DFO problems without modification. As noted in 
\cite{conn1998derivative}, and later corroborated in  \cite{shi2022numerical}, \texttt{IBO} methods are surprisingly robust in the presence of noise, for reasons that are not well understood. They achieve this without knowledge of the noise level in the function. Counter-examples exist showing that a standard IBO method can fail due to the effects of noise \cite{sun2022trust}, but such examples are rare.

%\sout{A potential strategy to mitigate the noise is to relax the ratio test \eqref{eq:ratio test} as in \cite{sun2022trust}. We leave this as a future research topic. }

\medskip\noindent{\it Feasible Interpolation Sets}. In cases when the function can only be evaluated inside the feasible region, all points in the initial interpolation set $Y_0$ must be feasible, something that can be expensive to compute. If the points in this set are not sufficiently well poised, model accuracy may be difficult or impossible to achieve, as discussed in  \cite{larson2019derivative}. We do not discuss here how to address this problem.

\medskip
 Considering the aforementioned points, the \texttt{FIBO} approach is well-suited for problems where the cost of evaluating the objective function $f$ outweighs that of solving the subproblem \eqref{outer}. The numerical experiments presented in section~\ref{sec:num}  support this observation.

\subsection{Related Work} \label{sec:literature}
Powell  played a pivotal role in the development of the \texttt{IBO} approach,  his contributions tracing back to \cite{Powe94b}. He  developed several methods and software for increasingly more complex problem classes.  
\texttt{NEWUOA} \cite{powell2006newuoa}, designed for unconstrained problems, employs the idea of minimum-Frobenius norm update to absorb degrees of freedom and reduce the number of interpolation points. \texttt{NEWUOA}  brings about significant reductions in linear algebra costs compared to previous implementations and represents a major contribution to the field.   \texttt{BOBYQA} \cite{powell2009bobyqa} extends this approach to bound constrained problems, while \texttt{LINCOA} \cite{powell2015fast} can handle linear constraints. These methods  are all accessible through the Python and Matlab interfaces developed by Ragonneau and Zhang \cite{PDFO}. 

Supporting theory for \texttt{IBO} methods as well as important algorithmic issues are discussed in  Conn et al. \cite{conn2009introduction}, which is a standard reference on derivative-free optimization. For a comprehensive review of recent developments in IBO methods see the recent survey by Larson et al. \cite{larson2019derivative}. 

\texttt{COBYLA} is an early method by Powell, designed for inequality constrained DFO problems  \cite{Powe94a}. It constructs linear models of the objective and constraints, and may be viewed as a predecessor of modern \texttt{IBO} methods. Ragonneau  \cite{ragonneau2022model} proposed an extension of \texttt{COBYLA}, named \texttt{COBYQA}, capable of handling  equality and inequality constraints.

The \texttt{FIBO} approach was first proposed by Conn et al. \cite{conn1998derivative}. In that study, the subproblem \eqref{outer} was solved using \texttt{NPSOL}, and the proposed method was compared against \texttt{COBYLA} and \texttt{LANCELOT}, the latter using  finite differences. However, the reported results are not correct, as \texttt{LANCELOT} computed finite difference estimates of the constraints as well, although these derivatives were available analytically. These unnecessary constraint evaluations were included in the total count, making it impossible to gleen the efficiency of \texttt{FIBO} compared with \texttt{LANCELOT}. Conejo et al. \cite{conejo2015trust} developed a \texttt{FIBO} method in which \eqref{outer} was solved using \texttt{ALGENCAN} \cite{doi:10.1137/060654797}. The numerical results are unsatisfactory too, as they tested against an inefficient method \cite{plantenga2009hopspack} based on a heuristic Augmented Lagrangian approach with subproblems solved using a pattern search method. \texttt{CONDOR} \cite{berghen2005condor} implements a \texttt{FIBO} method in which \eqref{outer} is solved using a specially designed SQP method but provides no numerical comparisons.

In terms of convergence analysis,  \cite{conejo2013global} establish convergence results for problems with convex constraints, assuming that the local model is always fully linear.  \cite{hough2022model} extend the global convergence analysis by proposing a generalized fully-linear model in the general convex constrained case, relaxing this fully-linear model assumption.  \cite{hough2022model} only provides numerical results on least-squares problems and not on constrained problems. 
%A similar idea as in the aformentioned methods is explored in \texttt{CONDOR}, which is designed to solve general inequality constraints by treating bound and linear constraints via active-set methods and the remaining with SQP\cite{berghen2005condor}.
%\jn{[What?]}

 %In summary, no convergence proof has yet been established for \texttt{FIBO} applied to problems with general nonlinear constraints, a topic we discuss next.

%Lewis et al. \cite{lewis2002globally} and Griffin et al. \cite{griffin2010nonlinearly} extend pattern search methods to handle nonlinear constraints by solving a sequence of unconstrained or linearly constrained problems, in which the objective includes a penalty or augmented Lagrangian reformulation with respect to the general nonlinear constraints . The subproblem is solved via pattern search methods but the derivative information of the constraints are not exploited. This approach does not yield an effective method. 

\section{Numerical Experiments}\label{sec:num}
We developed a Python implementation of the FIBO method, built around the unconstrained solver, \texttt{DFOTR}, maintained by Prof. Scheinberg's group and available on GitHub. The \texttt{DFOTR}  version used in our experiments implements a vanilla {IBO} method that does not include a geometry-improvement phase. Despite its simplicity,  \texttt{DFOTR} is robust and efficient, closely trailing in performance the state-of-the-art method, \texttt{NEWUOA}, in terms of function evaluations. To solve the trust region subproblem \eqref{outer} we employed \texttt{KNITRO}
\cite{ByrdNoceWalt06}.

There is no established constrained DFO solver that naturally lends itself for comparisons in our study.  \texttt{COBYLA} or \texttt{COBYQA} treat constraints as black boxes and model them by interpolation, putting them at a disadvantage in the setting of this paper. A more favorable alternative, as demonstrated in our experiments, is to utilize a standard code for nonlinear programming that approximates the gradient of the objective by finite differences and accepts analytic constraint derivatives.  The effectiveness of this approach to DFO is supported by the numerical study of Shi et al.
\cite{shi2022numerical} on inequality constrained DFO problems

The  details of the two solvers used in our experiments are as follows. 
\begin{itemize}
\item{\texttt{FIBO}.}
We employed  {\sc dfotr} as the base method, without altering its logic. The overall run was stopped if the trust region radius became smaller than $10^{-8}$ or if the number of objective function evaluations exceeded $500 \times \max(m,n)$, where $n$ is the number of variables and $m$ the number of constraints. We set the \texttt{stop\_predict} convergence test to 0 in {\sc dfotr} to avoid early termination caused by the ratio test. All remaining parameters in {\sc dfotr} were set to their default values. We use {\sc knitro}~12.4 with {exact gradients} to solve  \eqref{outer}, since exact gradient information is available for  this subproblem. We set \texttt{alg = 4} (SQP) and \texttt{hessopt = 6} {\sc (L-BFGS)} in  {\sc knitro} and set all other parameters to their default values. %\jn{[Why didn't we use {\sc knitro} with exact first and second derivatives? I think it would have given signficantly more efficient]}\mx{[I wanted to set hessopt to be the same for both approaches. Is it really going to be more significantly efficient if exact second derivatives are applied?]}

\item{\texttt{FD}.} We applied {\sc knitro}~12.4 directly to problem \eqref{problem}. The derivatives of the objective function $f$ were approximated by forward differences, wheras exact derivatives of the constraints were provided. We chose the following parameters in {\sc knitro:}
 \texttt{alg = 4} (SQP); \texttt{hessopt = 6} (L-BFGS) with memory size of 10;  \texttt{xtol = 1e-8, xtol\_iters = 1}, and \texttt{findiff\_terminate=0} (to obtain a similar convergence criterion as in \texttt{FIBO}). The maximum number of function evaluations allowed is, as before, $500 \times \max(m,n)$. 

\end{itemize}
We performed tests on a set of 38 general constrained problems  from the \texttt{CUTEst} collection.  Since the initial point provided by \texttt{CUTEst} is usually infeasible, we ran \texttt{KNITRO} to obtain a feasible point by replacing the objective by 0 in \eqref{problem}. The computational cost of obtaining this feasible point is not accounted for in the data presented below. 
%\jn{[Note that the observations remain the same when we include the cost of obtaining a feasible starting point for \texttt{FIBO} and run \texttt{FD} from an infeasible initial point $x_0$. See \ref{app: outerdfo infeas} for detailed numerical results.]} 
A fine implementation detail is that the initial iterate in {\sc dfotr} is the previously mentioned initial feasible point, rather than the interpolation point on the stencil with the minimum objective value. The reason for this choice is that such a point might be infeasible.
% \jn{[Is it that important to start with a feasible point? What is the issue otherwise? If we set the initial trust region radius to be quite large, the first run of the nonlinear solver should take care of thing.]}\mx{[If the initial trust region is set to be sufficiently large, the initial point is not required to be feasible. But I guess in practice, we do not know how big the initial trust region should be. Furthermore, I have noticed in my experiments that for the unconstrained DFO setting, choosing an initial region that is too large can lead to a less efficient algorithm in terms of function evaluations (the initial models are less accurate). ]}

\begin{figure}[htp]
	\centering
	\includegraphics[width = 0.35\textwidth]{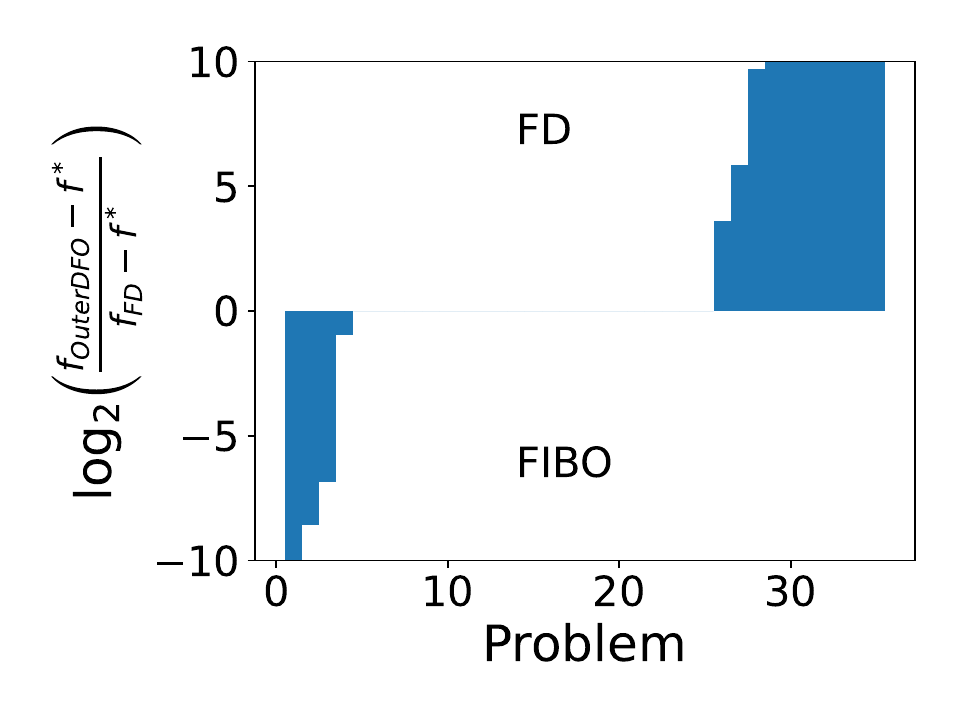}
	\caption{Log-ratio plot \eqref{acc} comparing the final accuracy  achieved by \texttt{FIBO} and \texttt{FD}. }
	\label{fig:acc}
\end{figure}

  In our first experiment, we assess the capability of \texttt{FIBO} to compute an accurate solution. We do so by recording the log-ratio  \cite{Mora02} given by
\begin{equation}\label{acc}
	\log_2 \left( \frac{\max(f_{\text{FIBO}}-f^*,10^{-8})}{\max(f_{\text{FD}} - f^*, 10^{-8})} \right),
\end{equation}
where $f_{\text{FIBO}}$ and $f_{\text{FD}}$ are the final objective values achieved by each method, up to 8 digits, and $f^*$ is the optimal value obtained by running \texttt{KNITRO} with exact gradients.
The results are depicted in Figure \ref{fig:acc}, showing the ratios \eqref{acc} plotted in ascending order. By design, the log-ratio ensures that a larger shaded region indicates a more successful method.
 We observe that \texttt{FD} slightly outperforms \texttt{FIBO} in terms of accuracy, which can be attributed to the fact that \texttt{FD} employs an approximation to the gradient and updates a quasi-Newton model, whereas \texttt{FIBO} constructs an interpolatory model of the objective function. For detailed numerical results,  refer to Appendix~\ref{app: outerdfo infeas}.

\begin{figure}[htp!]
	\centering
	\includegraphics[width=0.35\textwidth]{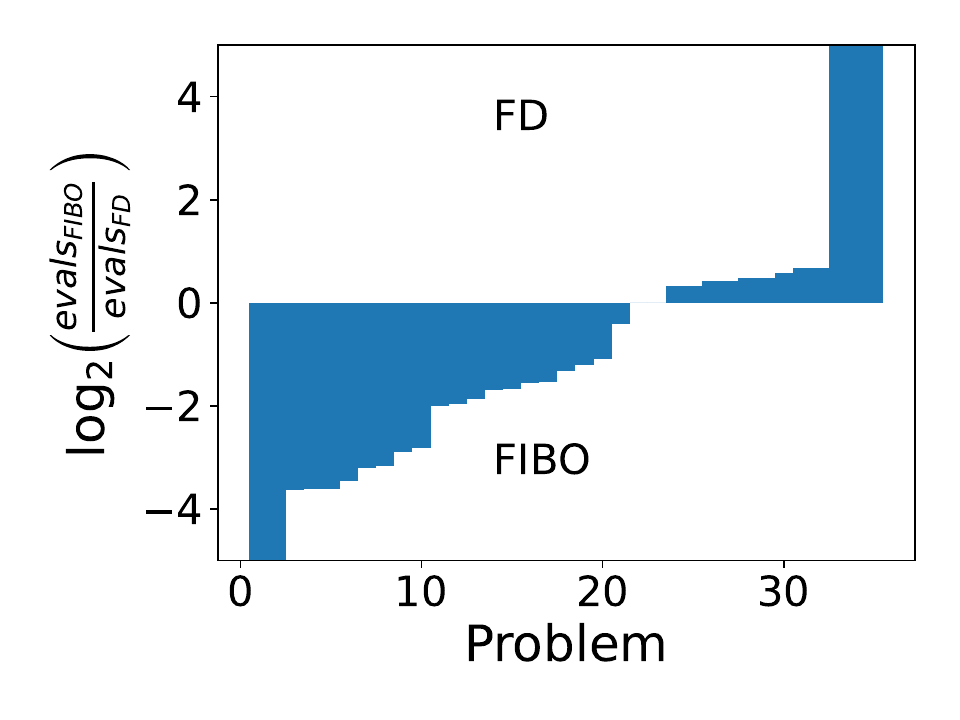}
	\includegraphics[width=0.35\textwidth]{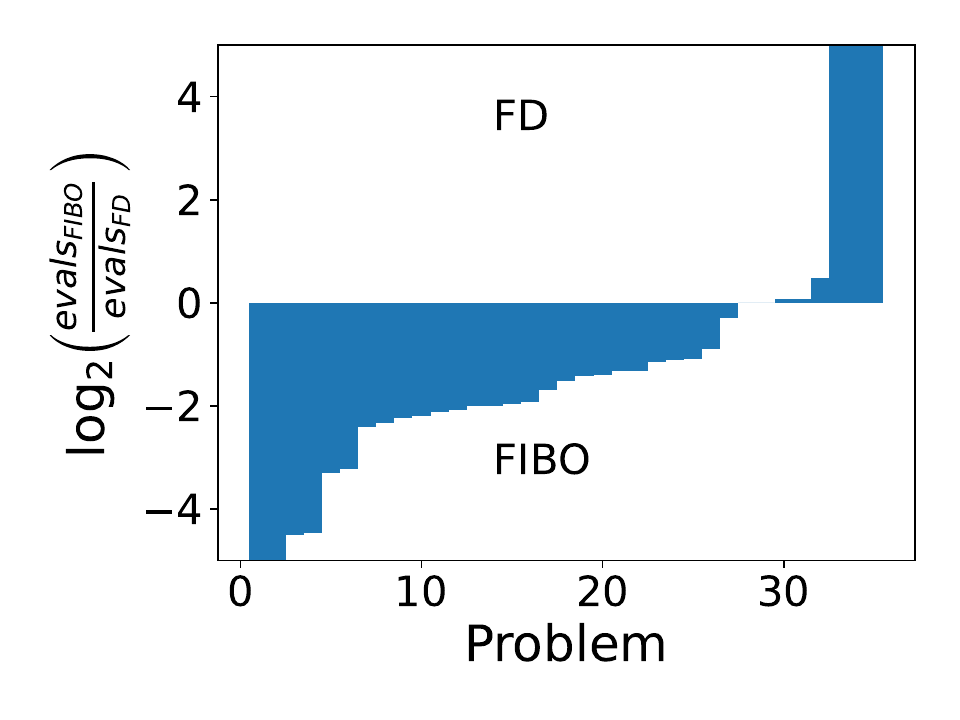} \\
	\includegraphics[width=0.35\textwidth]{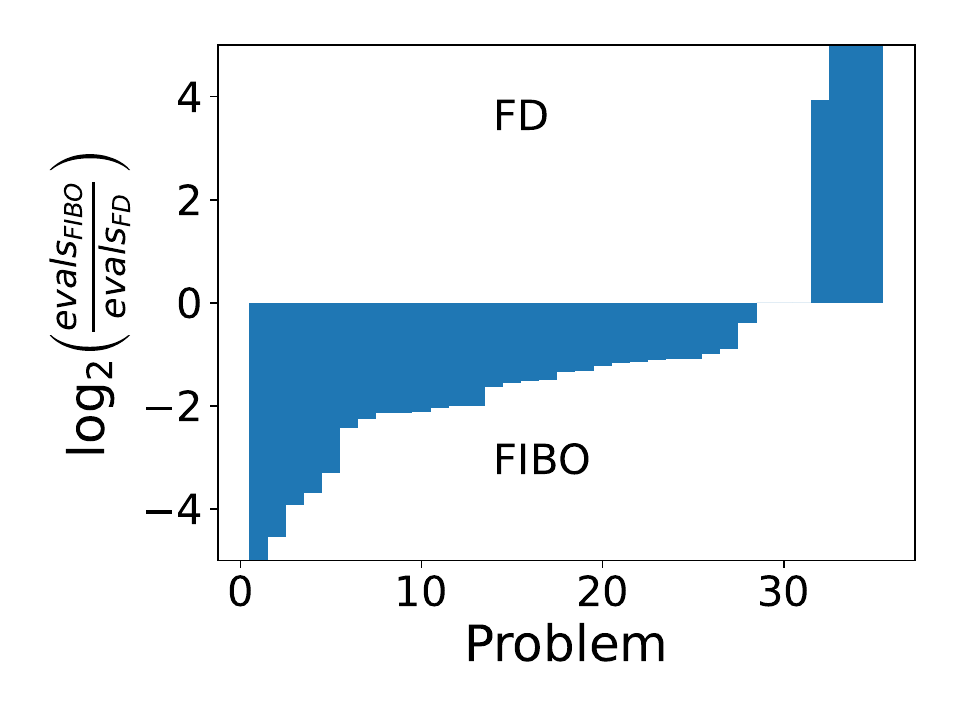}
	\includegraphics[width=0.35\textwidth]{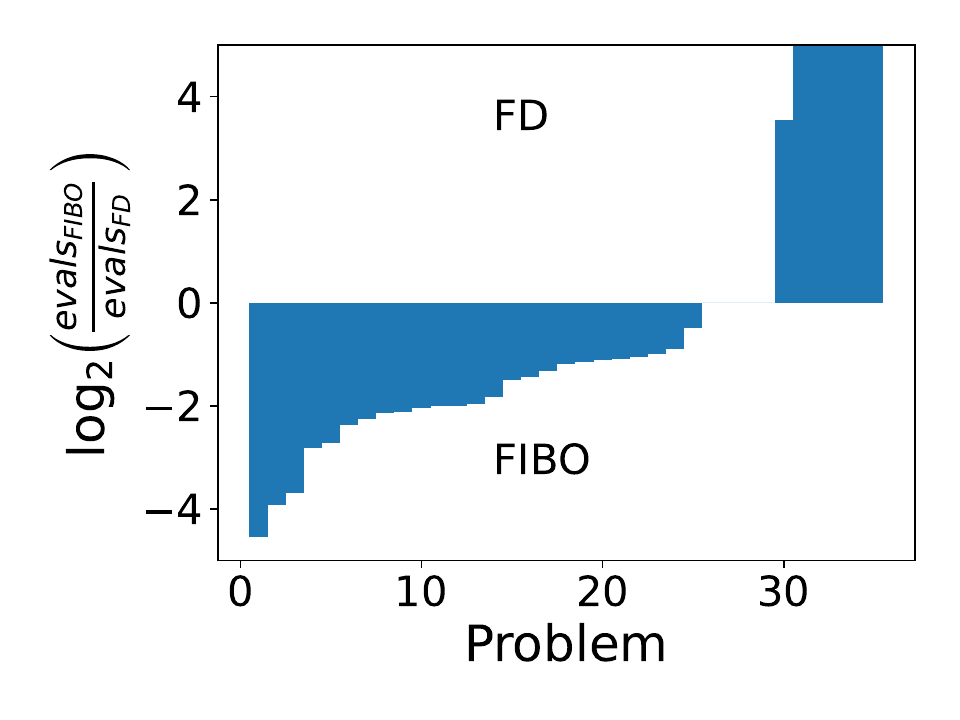}
	\caption{Log-ratio plot comparing \texttt{FIBO} and \texttt{FD} in terms of the number of objective function evaluations required to satisfy \eqref{efficiency} for 
%	$\tau = 10^{-1} \text{ (upper left)}, 10^{-3}\text{ (upper right)}, 10^{-5}\text{ (bottom left)}, 10^{-7}\text{ (bottom right)}$.}
	$\tau = 10^{-1}\text{ (upper left)},$ $ 10^{-3}\text{ (upper right)},\\$ $10^{-5}$ $\text{ (bottom left)}, 10^{-7}\text{ (bottom right)}$.}

	\label{feas plot f no failure}
\end{figure}

\begin{figure}[htp!]
	\centering
	\includegraphics[width=0.35\textwidth]{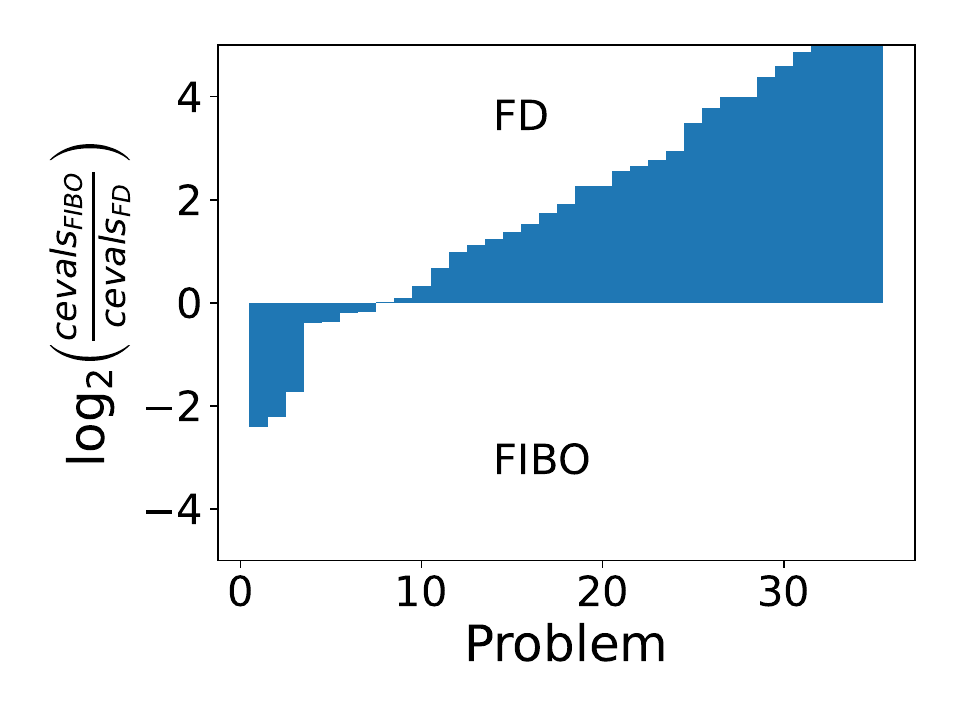}
	\includegraphics[width=0.35\textwidth]{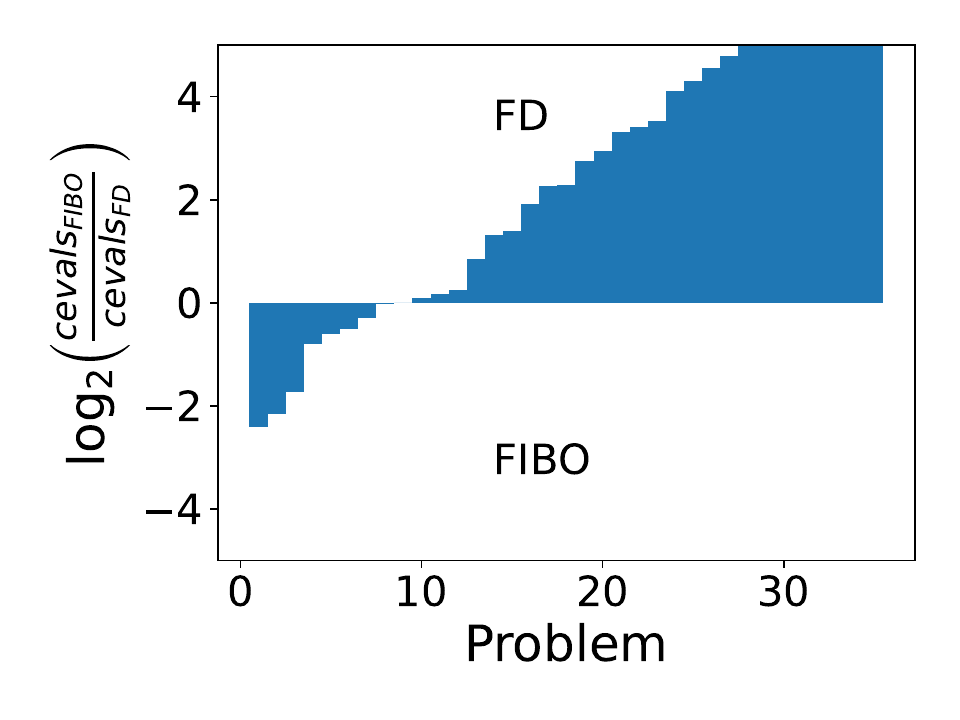} \\
	\includegraphics[width=0.35\textwidth]{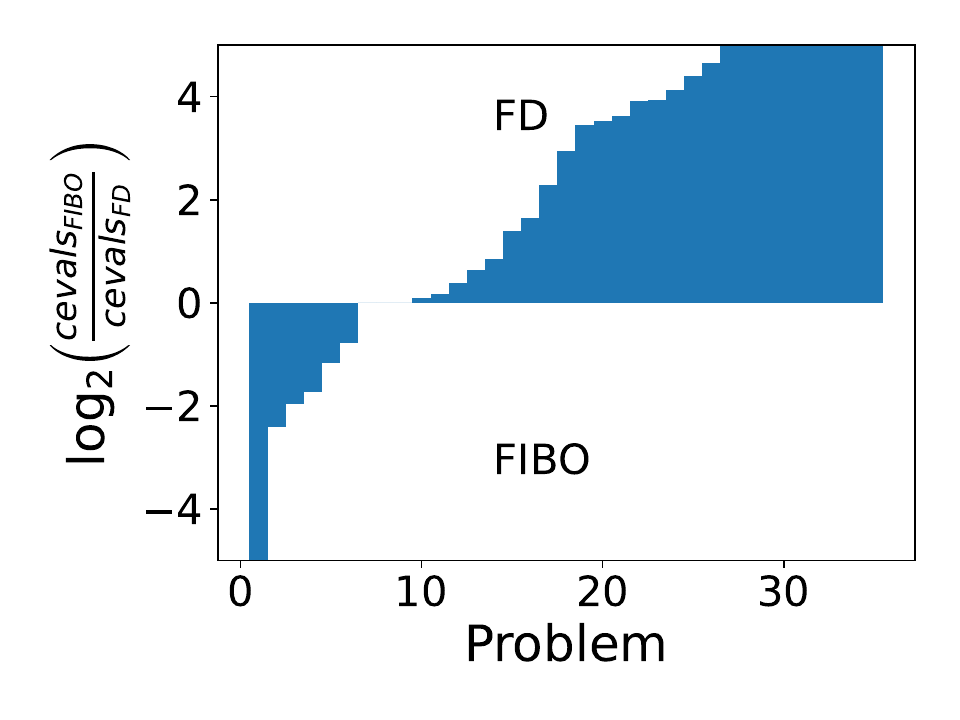}
	\includegraphics[width=0.35\textwidth]{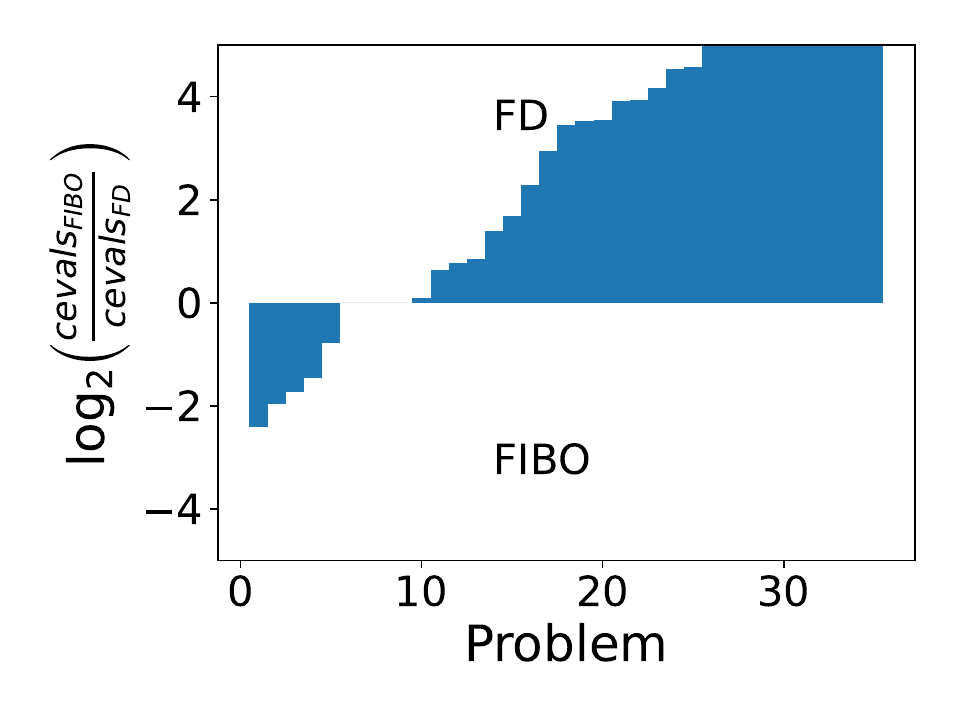}
	\caption{Log-ratio plot comparing \texttt{FIBO} and \texttt{FD} in terms of the number of constraint evaluations to satisfy \eqref{efficiency} for $\tau = 10^{-1} \text{ (upper left)}$, $10^{-3}\text{ (upper right)},$ $ 10^{-5}\text{ (bottom left)}$, $10^{-7}\text{ (bottom right)}$.}
	\label{feas plot c no failure}
\end{figure}

In the next set of experiments, we measure the efficiency of \texttt{FIBO} by recording the number of objective and constraint evaluations  required to achieve the accuracy
\begin{equation}\label{efficiency}
	f(x_k) \leq f^* + \tau \cdot \max(1, |f^*|),
\end{equation}
for $\tau = 10^{-1}, 10^{-3}, 10^{-5}$ and $10^{-7}$. We denote by $\text{evals}_{\text{FIBO}}$ and $\text{evals}_{\text{FD}}$ the number of objective function evaluations required by each method to satisfy \eqref{efficiency}. Similarly,  $\text{cevals}_{\text{FIBO}}$ and $\text{cevals}_{\text{FD}}$ denote the number of constraint evaluations. If  condition \eqref{efficiency} cannot be satisfied by an algorithm, we  set the corresponding quantity as a very large number.  We summarize the performance of  \texttt{IBO} and \texttt{FD}  in Figures~\ref{feas plot f no failure} and \ref{feas plot c no failure}, respectively, by reporting the log-ratios 
\begin{equation}
	\log_2 \left(\frac{\text{evals}_{\text{FIBO}}}{\text{evals}_{\text{FD}}}\right), \quad \log_2\left(\frac{\text{cevals}_{\text{FIBO}}}{\text{cevals}_{\text{FD}}}\right).
\end{equation} 

As shown in Figure \ref{feas plot f no failure}, \texttt{FIBO} outperforms \texttt{FD} in terms of objective evaluations across all accuracy levels, while \texttt{FD} is far more efficient in terms of constraint evaluations, as expected. This provides  support for the claim that \texttt{FIBO} is well-suited for problems in which the cost of the objective function dominates all other costs.

\section{Final Remarks}\label{sec:final}
In this paper, we studied the numerical efficiency  of the feasible interpolation-based optimization (\texttt{FIBO}) approach first proposed by Conn, Scheinberg and Toint \cite{conn1998derivative} for constrained optimization problems in which the gradient of the objective function is not available. While this method has been examined in previous studies \cite{conn1998derivative,conejo2015trust,berghen2005condor}, none of them  provided compelling numerical evidence regarding its practical viability. To the best of our knowledge, this method is rarely, if ever, utilized in practical applications, a situation this article attempts to rectify.
Our numerical tests  indicate that \texttt{FIBO} is a very competitive method for solving problems where the evaluation of the objective function $f$ dominates the cost of solving the trust region problem \eqref{outer}. This implies, in particular, that the evaluations of the constraints and their derivatives must be very inexpensive compared to the evaluation of $f$. This is the case when $f$ is the result of an expensive simulation while the constraints are provided analytically. Note, however, that \texttt{FIBO} does not require that the derivatives of the constraints be provided analytically; they could be approximated by finite differences. The only requirement is that these computations be inexpensive compared to an evaluation of $f$. 
  
%One attractive aspect of  \texttt{FIBO}  is that it is built by embedding a nonlinear optimization solver into an interpolation-based optimization (\texttt{IBO}) method for unconstrained optimization. The relative ease of doing so stands in constrast to the complexity associated with adapting \texttt{IBO} methods to handle constraints, as seen in  \texttt{LINCOA} \cite{powell2015fast}. Moreover, \texttt{FIBO} offers an added advantage in its ability to handle problems where the objective function is undefined outside of the feasible region. 
% The main challenge is in constructing an initial feasible interpolation set. While we can envision several strategies to address this issue, an examination of these approaches is beyond the scope of this paper.

 The primary limitation of \texttt{FIBO} lies in the need to solve the nonlinear programming subproblem \eqref{outer} to optimality at each iteration. Although we never encountered failures in our experiments, it cannot be assured that an optimal solution of this subproblem is always found. %\sout{While convergence guarantees are lacking for any method tackling the general constrained problem \eqref{problem}, in the case of \texttt{FIBO}, the potential large number of subproblems to be solved makes this lack of robustness a more pressing concern.}

In this paper, we have not explored scenarios in which the objective function exhibits noise. Nonetheless, it is our expectation that the \texttt{FIBO} approach will be effective in such conditions too. This view is based in the empirically observed robustness of the  \texttt{IBO} strategy to noise, as reported by Shi et al. \cite{shi2022numerical}.

% \mx{We demonstrated effectiveness of the \texttt{FIBO} approach based on empirical investigations. When only analytical constraints are present, our study has revealed that, with minor modifications, existing unconstrained IBO methods can be extended to solve general constrained DFO problems efficiently, by redirecting the cost of objective evaluations to
% additional constraint evaluations. Our framework applies to any IBO methods and gradient-based nonlinear programming solvers. Yet more sophisticated algorithmic design is required when the constraints cannot be violated, or black-box constraints are present. It is non-trivial to compute interpolation points that are feasible, when model-improving points are necessary to improve model quality. In addition, when constraints are expensive to evaluate, further algorithmic development for IBO is required to improve efficiency of the proposed method. We have observed fair performance of \texttt{FIBO} even if the termination criteria of the trust region subproblems \eqref{outer} are relaxed to a larger value. Yet it requires more careful investigations as to whether this is valid and stable. We hope our studies provide insights into the design of constrained DFO algorithms.} 

% However, we observed in our experiments that using only feasible points in interpolation did not produce singular interpolation linear system in most of the cases.
\bibliographystyle{plain}

\bibliography{../../References/references}
\appendix

\section{Numerical Results for Feasible Initial Points}\label{app:feas num}
\begin{landscape}
\begin{table}[htp]
\centering
\resizebox{0.77\columnwidth}{!}{\begin{tabular}{lll | llllllll | lllll}
\hline
\multicolumn{3}{l}{} & \multicolumn{8}{c}{FIBO}       & \multicolumn{5}{c}{knitro}         \\ \hline
Problem   & n   & m  & \#iter & \#feval & time & $f$ & feas err & \#iter(sub) & \#ceval & time(sub) &  \#iter & \#feval & time & $f$ & feas err \\ \hline
HS13 & $2$ & $1$ & $3$ & $5$  & $0.29$ & \num{ 1.000e+00} & \num{ 1.309e-14}&136&191&0.28& $20$ & $104$  & $0.09$ & \num{ 9.983e-01} & \num{ 6.579e-10} \\
HS22 & $2$ & $2$ & $2$ & $4$  & $0.04$ & \num{ 1.000e+00} & \num{ 0.000e+00}&10&25&0.03& $1$ & $8$  & $0.01$ & \num{ 1.000e+00} & \num{ 0.000e+00} \\
HS23 & $2$ & $5$ & $4$ & $6$  & $0.04$ & \num{ 2.000e+00} & \num{ 0.000e+00}&12&16&0.03& $6$ & $28$  & $0.03$ & \num{ 2.000e+00} & \num{ 0.000e+00} \\
HS26 & $3$ & $1$ & $53$ & $56$  & $12.96$ & \num{ 7.103e-10} & \num{ 2.220e-16}&5245&39137&12.85& $34$ & $244$  & $0.1$ & \num{ 7.434e-12} & \num{ 6.661e-16} \\
HS32 & $3$ & $2$ & $4$ & $7$  & $0.04$ & \num{ 1.000e+00} & \num{ 0.000e+00}&10&16&0.03& $2$ & $15$  & $0.01$ & \num{ 1.000e+00} & \num{ 0.000e+00} \\
HS34 & $3$ & $2$ & $6$ & $9$  & $0.1$ & \num{-8.340e-01} & \num{ 1.776e-15}&32&70&0.08& $7$ & $44$  & $0.02$ & \num{-8.340e-01} & \num{ 0.000e+00} \\
HS40 & $4$ & $3$ & $22$ & $26$  & $0.36$ & \num{-2.500e-01} & \num{ 1.110e-16}&56&225&0.26& $5$ & $44$  & $0.01$ & \num{-2.500e-01} & \num{ 1.110e-16} \\
HS44 & $4$ & $6$ & $3$ & $8$  & $0.04$ & \num{-1.000e+00}(*) & \num{ 0.000e+00}&12&23&0.03& $5$ & $36$  & $0.03$ & \num{-1.500e+01} & \num{ 0.000e+00} \\
HS47 & $5$ & $3$ & $23$ & $28$  & $3.67$ & \num{ 1.989e-07} & \num{ 3.781e-11}&1391&7878&3.47& $72$ & $696$  & $0.41$ & \num{ 1.007e+01} & \num{ 6.661e-16} \\
HS50 & $5$ & $3$ & $19$ & $24$  & $1.17$ & \num{ 1.056e+00} & \num{ 2.398e-14}&126&423&1.1& $16$ & $121$  & $0.04$ & \num{ 7.669e-17} & \num{ 4.441e-16} \\
HS64 & $3$ & $1$ & $43$ & $46$  & $1.04$ & \num{ 6.300e+03} & \num{ 0.000e+00}&300&662&0.95& $12$ & $65$  & $0.04$ & \num{ 6.300e+03} & \num{ 0.000e+00} \\
HS66 & $3$ & $2$ & $2$ & $6$  & $0.03$ & \num{ 5.182e-01}(*) & \num{ 3.563e-12}&7&12&0.02& $7$ & $43$  & $0.03$ & \num{ 5.182e-01} & \num{ 0.000e+00} \\
HS67 & $3$ & $14$ & $45$ & $49$  & $1.7$ & \num{-1.162e+03}(*) & \num{ 0.000e+00}&472&1408&1.59& $14$ & $75$  & $0.03$ & \num{-1.162e+03} & \num{ 0.000e+00} \\
HS72 & $4$ & $2$ & $7$ & $11$  & $0.26$ & \num{ 7.277e+02} & \num{ 2.819e-18}&44&99&0.21& $8$ & $55$  & $0.03$ & \num{ 7.277e+02} & \num{ 4.337e-19} \\
HS75 & $4$ & $5$ & $16$ & $20$  & $0.49$ & \num{ 5.174e+03} & \num{ 5.684e-14}&97&180&0.43& $6$ & $43$  & $0.04$ & \num{ 5.174e+03} & \num{ 8.527e-14} \\
HS85 & $5$ & $21$ & $31$ & $36$  & $1.46$ & \num{-2.216e+00} & \num{ 0.000e+00}&357&903&1.31& $11$ & $85$  & $0.06$ & \num{-2.216e+00} & \num{ 1.421e-14} \\
HS87 & $6$ & $4$ & $9$ & $15$  & $0.74$ & \num{ 8.997e+03} & \num{ 2.274e-13}&80&287&0.71& $14$ & $223$  & $0.07$ & \num{ 8.997e+03} & \num{ 2.416e-13} \\
HS88 & $2$ & $1$ & $5$ & $7$  & $0.4$ & \num{ 1.363e+00} & \num{ 0.000e+00}&92&318&0.39& $27$ & $156$  & $0.08$ & \num{ 1.363e+00} & \num{ 2.893e-17} \\
HS89 & $3$ & $1$ & $10$ & $13$  & $1.42$ & \num{ 1.363e+00} & \num{ 0.000e+00}&246&889&1.39& $5$ & $35$  & $0.04$ & \num{ 7.283e-08}$(\dagger)$ & \num{ 1.332e-01} \\
HS90 & $4$ & $1$ & $12$ & $16$  & $4.22$ & \num{ 1.363e+00} & \num{ 0.000e+00}&738&3531&4.18& $33$ & $285$  & $0.15$ & \num{ 1.363e+00} & \num{ 6.600e-18} \\
HS93 & $6$ & $2$ & $42$ & $48$  & $2.9$ & \num{ 1.351e+02} & \num{ 4.800e-12}&517&2433&2.61& $26$ & $295$  & $0.13$ & \num{ 1.351e+02} & \num{ 0.000e+00} \\
HS98 & $6$ & $4$ & $3$ & $10$  & $0.05$ & \num{ 3.136e+00}(*) & \num{ 1.776e-15}&5&8&0.02& $4$ & $40$  & $0.02$ & \num{ 3.136e+00} & \num{ 1.776e-15} \\
HS100 & $7$ & $4$ & $106$ & $113$  & $29.46$ & \num{ 6.806e+02} & \num{ 0.000e+00}&3280&22448&28.72& $51$ & $587$  & $0.21$ & \num{ 6.806e+02} & \num{ 7.161e-15} \\
HS101 & $7$ & $5$ & $103$ & $110$  & $19.42$ & \num{ 1.810e+03} & \num{ 3.955e-16}&2616&15697&18.39& $35$ & $354$  & $0.19$ & \num{ 1.810e+03} & \num{ 0.000e+00} \\
HS102 & $7$ & $5$ & $103$ & $110$  & $15.11$ & \num{ 9.119e+02} & \num{ 3.816e-16}&1985&10240&13.86& $21$ & $254$  & $0.13$ & \num{ 9.119e+02} & \num{ 2.498e-16} \\
HS103 & $7$ & $5$ & $56$ & $63$  & $5.5$ & \num{ 5.437e+02} & \num{ 1.318e-16}&825&4072&4.59& $20$ & $238$  & $0.1$ & \num{ 5.437e+02} & \num{ 2.359e-16} \\
HS104 & $8$ & $5$ & $245$ & $253$  & $43.5$ & \num{ 3.951e+00} & \num{ 9.714e-17}&6813&37268&41.25& $30$ & $359$  & $0.1$ & \num{ 3.951e+00} & \num{ 2.082e-16} \\
LOADBAL & $31$ & $31$ & $1$ & $15499$  & $0.93$ & \num{ 1.547e+00}(***) & \num{ 1.798e+308}&0&1&0.0& $40$ & $1380$  & $0.21$ & \num{ 4.529e-01} & \num{ 2.665e-15} \\
OPTPRLOC & $30$ & $30$ & $9$ & $39$  & $1.0$ & \num{-1.642e+01} & \num{ 1.990e-13}&54&127&0.45& $16$ & $583$  & $0.13$ & \num{-1.642e+01} & \num{ 1.421e-13} \\
CB3 & $3$ & $3$ & $1$ & $4$  & $0.01$ & \num{ 2.000e+00} & \num{ 0.000e+00}&2&3&0.0& $2$ & $15$  & $0.01$ & \num{ 2.000e+00} & \num{ 0.000e+00} \\
CRESC50 & $6$ & $100$ & $1000$ & $1006$  & $215.01$ & \num{ 5.939e-01}(**) & \num{ 5.482e-13}&26160&96691&210.07& $438$ & $5787$  & $1.76$ & \num{ 5.946e-01} & \num{ 1.188e-14} \\
DEMBO7 & $16$ & $20$ & $19$ & $35$  & $2.37$ & \num{ 1.748e+02} & \num{ 4.073e-15}&302&795&2.07& $17$ & $354$  & $0.16$ & \num{ 1.749e+02} & \num{ 6.661e-16} \\
DNIEPER & $57$ & $24$ & $24$ & $81$  & $7.72$ & \num{ 1.874e+04} & \num{ 5.684e-14}&273&1094&3.25& $8$ & $534$  & $0.11$ & \num{ 1.874e+04} & \num{ 5.684e-14} \\
EXPFITA & $5$ & $22$ & $1$ & $10999$  & $0.6$ & \num{ 2.999e+01}(***) & \num{ 1.798e+308}&0&1&0.0& $19$ & $140$  & $0.05$ & \num{ 1.137e-03} & \num{ 1.421e-14} \\
HIMMELBI & $100$ & $12$ & $114$ & $215$  & $72.23$ & \num{-1.467e+03}(*) & \num{ 1.789e-09}&2618&12586&67.95& $72$ & $7487$  & $0.72$ & \num{-1.736e+03} & \num{ 1.066e-14} \\
SYNTHES1 & $6$ & $6$ & $1$ & $2999$  & $0.21$ & \num{ 1.000e+01}(***) & \num{ 1.798e+308}&0&1&0.0& $9$ & $81$  & $0.03$ & \num{ 7.593e-01} & \num{ 1.110e-16} \\
TWOBARS & $2$ & $2$ & $7$ & $9$  & $0.1$ & \num{ 1.509e+00} & \num{ 0.000e+00}&24&108&0.09& $9$ & $57$  & $0.03$ & \num{ 1.509e+00} & \num{ 4.441e-16} \\
DIPIGRI & $7$ & $4$ & $116$ & $123$  & $85.85$ & \num{ 6.806e+02} & \num{ 3.581e-11}&7461&26264&84.59& $48$ & $581$  & $0.18$ & \num{ 6.806e+02} & \num{ 0.000e+00} \\
\end{tabular}}
\caption{Noiseless Problems with Feasible $x_0$;  $n$: number of variables, $m$: number of constraints, \#iter: number of (outer) iterations, \#feval: number of function evaluations, time: total CPU time passed,  $f$: final objective value, feas err: final feasibility error, \#iter(sub): total number of iterations for solving TR subproblem, \#ceval: total number of constraint evaluations(note that the number of constraint evaluations and function evaluations are same for \texttt{KNITRO}), time(sub): total CPU time elapsed for solving subproblem. * indicates that \texttt{FIBO} terminates with singular interpolation system error, ** indicates that \texttt{FIBO} terminates with maximum number of iterations, *** indicates that \texttt{FIBO} terminates with maximum number of function evaluations. }
\label{tab:noiseless acc feas}
\end{table}

\end{landscape}

\begin{landscape}
%tau = 1e-1
\begin{table}[htp]
\centering
\resizebox{0.77\columnwidth}{!}{\begin{tabular}{lll | llllllll | lllll}
\hline
\multicolumn{3}{l}{} & \multicolumn{8}{c}{FIBO}       & \multicolumn{5}{c}{knitro}         \\ \hline
Problem   & n   & m  & \#iter & \#feval & time & $f$ & feas err & \#iter(sub) & \#ceval & time(sub) &  \#iter & \#feval & time & $f$ & feas err \\ \hline
HS13 & $2$ & $1$ & $2$ & $4$  & $0.2$ & \num{ 1.000e+00} & \num{ 1.310e-14}&75&104&0.19& $8$ & $36$  & $0.02$ & \num{ 1.080e+00} & \num{ 0.000e+00} \\
HS22 & $2$ & $2$ & $1$ & $3$  & $0.03$ & \num{ 1.000e+00} & \num{ 1.789e-12}&9&23&0.02& $0$ & $3$  & $0.0$ & \num{ 1.000e+00} & \num{ 1.843e-14} \\
HS23 & $2$ & $5$ & $3$ & $5$  & $0.05$ & \num{ 2.000e+00} & \num{ 0.000e+00}&11&14&0.03& $3$ & $16$  & $0.01$ & \num{ 2.040e+00} & \num{ 0.000e+00} \\
HS26 & $3$ & $1$ & $4$ & $7$  & $0.09$ & \num{ 1.996e-02} & \num{ 9.592e-14}&32&66&0.07& $14$ & $87$  & $0.04$ & \num{ 4.252e-07} & \num{ 9.633e-07} \\
HS32 & $3$ & $2$ & $4$ & $7$  & $0.09$ & \num{ 1.000e+00} & \num{ 0.000e+00}&10&16&0.04& $2$ & $15$  & $0.02$ & \num{ 1.000e+00} & \num{ 0.000e+00} \\
HS34 & $3$ & $2$ & $5$ & $8$  & $0.09$ & \num{-7.820e-01} & \num{ 1.066e-14}&31&68&0.08& $4$ & $29$  & $0.02$ & \num{-8.266e-01} & \num{ 0.000e+00} \\
HS40 & $4$ & $3$ & $1$ & $5$  & $0.02$ & \num{-2.454e-01} & \num{ 2.554e-15}&5&15&0.01& $0$ & $4$  & $0.0$ & \num{-2.454e-01} & \num{ 7.334e-09} \\
HS44 & $4$ & $6$ & $3$ & $8$  & $0.09$ & \num{-1.000e+00}(*) & \num{ 0.000e+00}&12&23&0.07& $4$ & $30$  & $0.02$ & \num{-1.500e+01} & \num{ 0.000e+00} \\
HS47 & $5$ & $3$ & $6$ & $11$  & $0.14$ & \num{ 5.639e-02} & \num{ 1.466e-12}&46&118&0.12& $17$ & $134$  & $0.04$ & \num{ 4.638e-08} & \num{ 3.167e-09} \\
HS50 & $5$ & $3$ & $19$ & $24$  & $1.26$ & \num{ 1.056e+00} & \num{ 2.398e-14}&126&423&1.19& $9$ & $71$  & $0.03$ & \num{ 4.941e-02} & \num{ 4.441e-16} \\
HS64 & $3$ & $1$ & $1$ & $4$  & $0.05$ & \num{ 6.813e+03} & \num{ 4.742e-11}&14&72&0.05& $0$ & $3$  & $0.0$ & \num{ 6.813e+03} & \num{ 0.000e+00} \\
HS66 & $3$ & $2$ & $1$ & $4$  & $0.02$ & \num{ 5.800e-01} & \num{ 3.563e-12}&6&10&0.01& $0$ & $3$  & $0.0$ & \num{ 5.800e-01} & \num{ 0.000e+00} \\
HS67 & $3$ & $14$ & $31$ & $34$  & $1.91$ & \num{-1.053e+03} & \num{ 0.000e+00}&402&1319&1.82& $8$ & $45$  & $0.02$ & \num{-1.157e+03} & \num{ 0.000e+00} \\
HS72 & $4$ & $2$ & $1$ & $5$  & $0.01$ & \num{ 7.407e+02} & \num{ 0.000e+00}&2&5&0.01& $0$ & $4$  & $0.0$ & \num{ 7.407e+02} & \num{ 0.000e+00} \\
HS75 & $4$ & $5$ & $1$ & $5$  & $0.02$ & \num{ 5.350e+03} & \num{ 1.259e-08}&4&8&0.01& $5$ & $37$  & $0.03$ & \num{ 5.174e+03} & \num{ 7.739e-08} \\
HS85 & $5$ & $21$ & $29$ & $34$  & $1.22$ & \num{-2.061e+00} & \num{ 1.421e-14}&348&880&1.04& $10$ & $78$  & $0.05$ & \num{-2.216e+00} & \num{ 2.389e-08} \\
HS87 & $6$ & $4$ & $1$ & $7$  & $0.07$ & \num{ 8.997e+03} & \num{ 2.274e-13}&9&24&0.06& $0$ & $5$  & $0.0$ & \num{ 8.997e+03} & \num{ 2.956e-12} \\
HS88 & $2$ & $1$ & $1$ & $3$  & $0.07$ & \num{ 1.385e+00} & \num{ 2.158e-14}&19&41&0.07& $0$ & $3$  & $0.0$ & \num{ 1.385e+00} & \num{ 0.000e+00} \\
HS89 & $3$ & $1$ & $2$ & $5$  & $0.29$ & \num{ 1.365e+00} & \num{ 1.123e-16}&57&205&0.28& $5$ & $35$  & $0.03$ & \num{ 7.283e-08}$(\dagger)$ & \num{ 1.332e-01} \\
HS90 & $4$ & $1$ & $2$ & $6$  & $0.53$ & \num{ 1.446e+00} & \num{ 6.661e-16}&86&454&0.52& $24$ & $228$  & $0.11$ & \num{ 1.362e+00} & \num{ 5.472e-07} \\
HS93 & $6$ & $2$ & $1$ & $7$  & $0.11$ & \num{ 1.371e+02} & \num{ 2.946e-11}&10&24&0.11& $0$ & $5$  & $0.0$ & \num{ 1.371e+02} & \num{ 0.000e+00} \\
HS98 & $6$ & $4$ & $2$ & $8$  & $0.03$ & \num{ 3.136e+00} & \num{ 0.000e+00}&4&6&0.02& $3$ & $32$  & $0.02$ & \num{ 3.136e+00} & \num{ 0.000e+00} \\
HS100 & $7$ & $4$ & $1$ & $8$  & $0.12$ & \num{ 7.140e+02} & \num{ 0.000e+00}&22&79&0.12& $0$ & $5$  & $0.0$ & \num{ 7.140e+02} & \num{ 0.000e+00} \\
HS101 & $7$ & $5$ & $26$ & $33$  & $5.24$ & \num{ 1.982e+03} & \num{ 3.053e-16}&740&4383&5.04& $9$ & $106$  & $0.04$ & \num{ 1.848e+03} & \num{ 0.000e+00} \\
HS102 & $7$ & $5$ & $34$ & $41$  & $4.55$ & \num{ 9.948e+02} & \num{ 2.776e-17}&641&3320&4.29& $11$ & $159$  & $0.05$ & \num{ 9.191e+02} & \num{ 0.000e+00} \\
HS103 & $7$ & $5$ & $15$ & $22$  & $1.69$ & \num{ 5.623e+02} & \num{ 8.913e-12}&276&1263&1.63& $16$ & $202$  & $0.06$ & \num{ 5.437e+02} & \num{ 1.040e-07} \\
HS104 & $8$ & $5$ & $1$ & $9$  & $0.11$ & \num{ 4.200e+00} & \num{ 1.459e-11}&14&41&0.11& $0$ & $6$  & $0.0$ & \num{ 4.200e+00} & \num{ 0.000e+00} \\
LOADBAL & $31$ & $31$ & $1$ & $15499$  & $0.97$ & \num{ 1.547e+00}(***) & \num{ 1.798e+308}&0&1&0.0& $6$ & $231$  & $0.04$ & \num{ 5.252e-01} & \num{ 1.710e-14} \\
OPTPRLOC & $30$ & $30$ & $5$ & $35$  & $0.46$ & \num{-1.567e+01} & \num{ 2.220e-16}&44&104&0.4& $2$ & $103$  & $0.04$ & \num{-1.587e+01} & \num{ 6.661e-16} \\
CB3 & $3$ & $3$ & $1$ & $4$  & $0.01$ & \num{ 2.000e+00} & \num{ 0.000e+00}&2&3&0.01& $1$ & $10$  & $0.01$ & \num{ 2.000e+00} & \num{ 7.270e-13} \\
CRESC50 & $6$ & $100$ & $16$ & $22$  & $2.44$ & \num{ 6.057e-01} & \num{ 9.032e-11}&344&1528&2.38& $58$ & $700$  & $0.32$ & \num{ 6.155e-01} & \num{ 0.000e+00} \\
DEMBO7 & $16$ & $20$ & $13$ & $29$  & $1.39$ & \num{ 1.879e+02} & \num{ 5.917e-08}&198&507&1.19& $15$ & $318$  & $0.1$ & \num{ 1.749e+02} & \num{ 1.869e-13} \\
DNIEPER & $57$ & $24$ & $4$ & $61$  & $1.82$ & \num{ 1.941e+04} & \num{ 4.264e-11}&96&457&1.35& $2$ & $177$  & $0.03$ & \num{ 1.874e+04} & \num{ 1.800e-07} \\
EXPFITA & $5$ & $22$ & $1$ & $10999$  & $0.61$ & \num{ 2.999e+01}(***) & \num{ 1.798e+308}&0&1&0.0& $7$ & $56$  & $0.02$ & \num{ 9.889e-02} & \num{ 1.354e-14} \\
HIMMELBI & $100$ & $12$ & $114$ & $215$  & $70.47$ & \num{-1.467e+03}(*) & \num{ 1.789e-09}&2618&12586&68.05& $3$ & $408$  & $0.09$ & \num{-1.643e+03} & \num{ 0.000e+00} \\
SYNTHES1 & $6$ & $6$ & $1$ & $2999$  & $0.22$ & \num{ 1.000e+01}(***) & \num{ 1.798e+308}&0&1&0.0& $7$ & $65$  & $0.02$ & \num{ 7.593e-01} & \num{ 3.261e-07} \\
TWOBARS & $2$ & $2$ & $2$ & $4$  & $0.05$ & \num{ 1.523e+00} & \num{ 7.885e-09}&13&38&0.04& $7$ & $49$  & $0.05$ & \num{ 1.509e+00} & \num{ 2.179e-09} \\
DIPIGRI & $7$ & $4$ & $1$ & $8$  & $0.14$ & \num{ 7.140e+02} & \num{ 0.000e+00}&22&79&0.13& $0$ & $5$  & $0.0$ & \num{ 7.140e+02} & \num{ 0.000e+00} \\
\end{tabular}}
\caption{Noiseless Problems with Feasible $x_0$. $\tau = 10^{-1}$; $n$: number of variables, $m$: number of constraints, \#iter: number of (outer) iterations, \#feval: number of function evaluations, time: total CPU time passed,  $f$: final objective value, feas err: final feasibility error, \#iter(sub): total number of iterations for solving TR subproblem, \#ceval: total number of constraint evaluations(note that the number of constraint evaluations and function evaluations are same for \texttt{KNITRO}), time(sub): total CPU time elapsed for solving subproblem. * indicates that \texttt{FIBO} terminates with singular interpolation system error, ** indicates that \texttt{FIBO} terminates with maximum number of iterations, *** indicates that \texttt{FIBO} terminates with maximum number of function evaluations. }
\label{tab:noiseless feas 1e-1}
\end{table}

%tau = 1e-3
\begin{table}[htp]
\centering
\resizebox{0.77\columnwidth}{!}{\begin{tabular}{lll | llllllll | lllll}
\hline
\multicolumn{3}{l}{} & \multicolumn{8}{c}{FIBO}       & \multicolumn{5}{c}{knitro}         \\ \hline
Problem   & n   & m  & \#iter & \#feval & time & $f$ & feas err & \#iter(sub) & \#ceval & time(sub) &  \#iter & \#feval & time & $f$ & feas err \\ \hline
HS13 & $2$ & $1$ & $2$ & $4$  & $0.18$ & \num{ 1.000e+00} & \num{ 1.310e-14}&75&104&0.17& $18$ & $88$  & $0.08$ & \num{ 1.001e+00} & \num{ 0.000e+00} \\
HS22 & $2$ & $2$ & $1$ & $3$  & $0.03$ & \num{ 1.000e+00} & \num{ 1.789e-12}&9&23&0.02& $0$ & $3$  & $0.0$ & \num{ 1.000e+00} & \num{ 1.843e-14} \\
HS23 & $2$ & $5$ & $3$ & $5$  & $0.04$ & \num{ 2.000e+00} & \num{ 0.000e+00}&11&14&0.03& $4$ & $20$  & $0.02$ & \num{ 2.000e+00} & \num{ 0.000e+00} \\
HS26 & $3$ & $1$ & $36$ & $39$  & $1.85$ & \num{ 1.737e-04} & \num{ 1.040e-12}&754&4338&1.77& $14$ & $87$  & $0.03$ & \num{ 4.252e-07} & \num{ 9.633e-07} \\
HS32 & $3$ & $2$ & $4$ & $7$  & $0.05$ & \num{ 1.000e+00} & \num{ 0.000e+00}&10&16&0.03& $2$ & $15$  & $0.02$ & \num{ 1.000e+00} & \num{ 0.000e+00} \\
HS34 & $3$ & $2$ & $6$ & $9$  & $0.11$ & \num{-8.340e-01} & \num{ 1.776e-15}&32&70&0.08& $6$ & $39$  & $0.01$ & \num{-8.340e-01} & \num{ 5.446e-09} \\
HS40 & $4$ & $3$ & $2$ & $6$  & $0.03$ & \num{-2.491e-01} & \num{ 2.189e-09}&8&21&0.02& $3$ & $32$  & $0.02$ & \num{-2.500e-01} & \num{ 2.040e-09} \\
HS44 & $4$ & $6$ & $3$ & $8$  & $0.05$ & \num{-1.000e+00}(*) & \num{ 0.000e+00}&12&23&0.03& $4$ & $30$  & $0.03$ & \num{-1.500e+01} & \num{ 0.000e+00} \\
HS47 & $5$ & $3$ & $23$ & $28$  & $3.41$ & \num{ 1.989e-07} & \num{ 3.781e-11}&1391&7878&3.33& $72$ & $696$  & $0.35$ & \num{ 1.007e+01} & \num{ 6.661e-16} \\
HS50 & $5$ & $3$ & $19$ & $24$  & $1.02$ & \num{ 1.056e+00} & \num{ 2.398e-14}&126&423&0.96& $11$ & $85$  & $0.03$ & \num{ 1.339e-06} & \num{ 8.882e-16} \\
HS64 & $3$ & $1$ & $13$ & $16$  & $0.46$ & \num{ 6.303e+03} & \num{ 1.443e-15}&126&399&0.44& $7$ & $40$  & $0.02$ & \num{ 6.300e+03} & \num{ 7.051e-10} \\
HS66 & $3$ & $2$ & $2$ & $5$  & $0.02$ & \num{ 5.182e-01} & \num{ 3.563e-12}&7&12&0.02& $3$ & $21$  & $0.02$ & \num{ 5.186e-01} & \num{ 0.000e+00} \\
HS67 & $3$ & $14$ & $38$ & $41$  & $1.51$ & \num{-1.162e+03} & \num{ 0.000e+00}&448&1377&1.42& $9$ & $50$  & $0.02$ & \num{-1.161e+03} & \num{ 0.000e+00} \\
HS72 & $4$ & $2$ & $5$ & $9$  & $0.17$ & \num{ 7.284e+02} & \num{ 0.000e+00}&38&90&0.16& $3$ & $24$  & $0.02$ & \num{ 7.277e+02} & \num{ 1.663e-07} \\
HS75 & $4$ & $5$ & $16$ & $20$  & $0.32$ & \num{ 5.174e+03} & \num{ 5.684e-14}&97&180&0.28& $5$ & $37$  & $0.04$ & \num{ 5.174e+03} & \num{ 7.739e-08} \\
HS85 & $5$ & $21$ & $31$ & $36$  & $0.95$ & \num{-2.216e+00} & \num{ 0.000e+00}&357&903&0.87& $10$ & $78$  & $0.05$ & \num{-2.216e+00} & \num{ 2.389e-08} \\
HS87 & $6$ & $4$ & $1$ & $7$  & $0.07$ & \num{ 8.997e+03} & \num{ 2.274e-13}&9&24&0.06& $0$ & $5$  & $0.0$ & \num{ 8.997e+03} & \num{ 2.956e-12} \\
HS88 & $2$ & $1$ & $2$ & $4$  & $0.2$ & \num{ 1.363e+00} & \num{ 9.216e-16}&41&139&0.19& $23$ & $140$  & $0.06$ & \num{ 1.362e+00} & \num{ 4.623e-07} \\
HS89 & $3$ & $1$ & $4$ & $7$  & $0.49$ & \num{ 1.363e+00} & \num{ 3.663e-16}&101&374&0.48& $5$ & $35$  & $0.04$ & \num{ 7.283e-08}$(\dagger)$ & \num{ 1.332e-01} \\
HS90 & $4$ & $1$ & $6$ & $10$  & $1.54$ & \num{ 1.363e+00} & \num{ 5.260e-17}&277&1541&1.52& $24$ & $228$  & $0.14$ & \num{ 1.362e+00} & \num{ 5.472e-07} \\
HS93 & $6$ & $2$ & $21$ & $27$  & $1.5$ & \num{ 1.352e+02} & \num{ 4.695e-11}&333&1706&1.43& $8$ & $87$  & $0.04$ & \num{ 1.351e+02} & \num{ 0.000e+00} \\
HS98 & $6$ & $4$ & $2$ & $8$  & $0.04$ & \num{ 3.136e+00} & \num{ 0.000e+00}&4&6&0.01& $3$ & $32$  & $0.02$ & \num{ 3.136e+00} & \num{ 0.000e+00} \\
HS100 & $7$ & $4$ & $51$ & $58$  & $15.18$ & \num{ 6.813e+02} & \num{ 1.388e-17}&1609&11182&15.01& $3$ & $55$  & $0.02$ & \num{ 6.811e+02} & \num{ 0.000e+00} \\
HS101 & $7$ & $5$ & $44$ & $51$  & $7.81$ & \num{ 1.812e+03} & \num{ 3.006e-10}&1239&7244&7.65& $18$ & $198$  & $0.1$ & \num{ 1.810e+03} & \num{ 2.139e-08} \\
HS102 & $7$ & $5$ & $69$ & $76$  & $10.13$ & \num{ 9.125e+02} & \num{ 4.272e-13}&1483&7663&9.89& $15$ & $200$  & $0.06$ & \num{ 9.119e+02} & \num{ 5.354e-07} \\
HS103 & $7$ & $5$ & $37$ & $44$  & $3.94$ & \num{ 5.437e+02} & \num{ 7.980e-16}&671&3465&3.8& $16$ & $202$  & $0.08$ & \num{ 5.437e+02} & \num{ 1.040e-07} \\
HS104 & $8$ & $5$ & $47$ & $55$  & $6.18$ & \num{ 3.952e+00} & \num{ 5.218e-15}&1210&6047&5.99& $20$ & $259$  & $0.05$ & \num{ 3.951e+00} & \num{ 1.537e-07} \\
LOADBAL & $31$ & $31$ & $1$ & $15499$  & $0.9$ & \num{ 1.547e+00}(***) & \num{ 1.798e+308}&0&1&0.0& $9$ & $330$  & $0.05$ & \num{ 4.531e-01} & \num{ 1.421e-14} \\
OPTPRLOC & $30$ & $30$ & $6$ & $36$  & $0.44$ & \num{-1.642e+01} & \num{ 0.000e+00}&47&111&0.39& $3$ & $136$  & $0.06$ & \num{-1.641e+01} & \num{ 0.000e+00} \\
CB3 & $3$ & $3$ & $1$ & $4$  & $0.01$ & \num{ 2.000e+00} & \num{ 0.000e+00}&2&3&0.01& $1$ & $10$  & $0.0$ & \num{ 2.000e+00} & \num{ 7.270e-13} \\
CRESC50 & $6$ & $100$ & $68$ & $74$  & $18.0$ & \num{ 5.942e-01} & \num{ 0.000e+00}&2416&11257&17.76& $438$ & $5787$  & $1.81$ & \num{ 5.946e-01} & \num{ 1.188e-14} \\
DEMBO7 & $16$ & $20$ & $18$ & $34$  & $2.11$ & \num{ 1.749e+02} & \num{ 1.067e-10}&301&793&2.03& $15$ & $318$  & $0.1$ & \num{ 1.749e+02} & \num{ 1.869e-13} \\
DNIEPER & $57$ & $24$ & $5$ & $62$  & $1.42$ & \num{ 1.874e+04} & \num{ 1.405e-09}&100&465&1.35& $2$ & $177$  & $0.06$ & \num{ 1.874e+04} & \num{ 1.800e-07} \\
EXPFITA & $5$ & $22$ & $1$ & $10999$  & $0.58$ & \num{ 2.999e+01}(***) & \num{ 1.798e+308}&0&1&0.0& $12$ & $91$  & $0.03$ & \num{ 1.355e-03} & \num{ 6.750e-14} \\
HIMMELBI & $100$ & $12$ & $114$ & $215$  & $68.58$ & \num{-1.467e+03}(*) & \num{ 1.789e-09}&2618&12586&67.1& $12$ & $1326$  & $0.26$ & \num{-1.734e+03} & \num{ 5.684e-14} \\
SYNTHES1 & $6$ & $6$ & $1$ & $2999$  & $0.21$ & \num{ 1.000e+01}(***) & \num{ 1.798e+308}&0&1&0.0& $7$ & $65$  & $0.02$ & \num{ 7.593e-01} & \num{ 3.261e-07} \\
TWOBARS & $2$ & $2$ & $3$ & $5$  & $0.06$ & \num{ 1.509e+00} & \num{ 4.742e-11}&16&55&0.05& $7$ & $49$  & $0.04$ & \num{ 1.509e+00} & \num{ 2.179e-09} \\
DIPIGRI & $7$ & $4$ & $51$ & $58$  & $11.33$ & \num{ 6.812e+02} & \num{ 6.939e-18}&1491&9509&11.16& $3$ & $55$  & $0.01$ & \num{ 6.811e+02} & \num{ 0.000e+00} \\
\end{tabular}}
\caption{Noiseless Problems with Feasible $x_0$. $\tau = 10^{-3}$; $n$: number of variables, $m$: number of constraints, \#iter: number of (outer) iterations, \#feval: number of function evaluations, time: total CPU time passed,  $f$: final objective value, feas err: final feasibility error, \#iter(sub): total number of iterations for solving TR subproblem, \#ceval: total number of constraint evaluations(note that the number of constraint evaluations and function evaluations are same for \texttt{KNITRO}), time(sub): total CPU time elapsed for solving subproblem. * indicates that \texttt{FIBO} terminates with singular interpolation system error, ** indicates that \texttt{FIBO} terminates with maximum number of iterations, *** indicates that \texttt{FIBO} terminates with maximum number of function evaluations. }
\label{tab:noiseless feas 1e-3}
\end{table}

%tau = 1e-5
\begin{table}[htp]
\centering
\resizebox{0.77\columnwidth}{!}{\begin{tabular}{lll | llllllll | lllll}
\hline
\multicolumn{3}{l}{} & \multicolumn{8}{c}{FIBO}       & \multicolumn{5}{c}{knitro}         \\ \hline
Problem   & n   & m  & \#iter & \#feval & time & $f$ & feas err & \#iter(sub) & \#ceval & time(sub) &  \#iter & \#feval & time & $f$ & feas err \\ \hline
HS13 & $2$ & $1$ & $3$ & $5$  & $0.29$ & \num{ 1.000e+00} & \num{ 1.309e-14}&136&191&0.28& $19$ & $100$  & $0.1$ & \num{ 9.983e-01} & \num{ 6.579e-10} \\
HS22 & $2$ & $2$ & $1$ & $3$  & $0.03$ & \num{ 1.000e+00} & \num{ 1.789e-12}&9&23&0.02& $0$ & $3$  & $0.0$ & \num{ 1.000e+00} & \num{ 1.843e-14} \\
HS23 & $2$ & $5$ & $3$ & $5$  & $0.03$ & \num{ 2.000e+00} & \num{ 0.000e+00}&11&14&0.02& $5$ & $24$  & $0.02$ & \num{ 2.000e+00} & \num{ 0.000e+00} \\
HS26 & $3$ & $1$ & $38$ & $41$  & $8.14$ & \num{ 8.872e-07} & \num{ 4.441e-16}&3356&26976&8.07& $14$ & $87$  & $0.03$ & \num{ 4.252e-07} & \num{ 9.633e-07} \\
HS32 & $3$ & $2$ & $4$ & $7$  & $0.04$ & \num{ 1.000e+00} & \num{ 0.000e+00}&10&16&0.03& $2$ & $15$  & $0.02$ & \num{ 1.000e+00} & \num{ 0.000e+00} \\
HS34 & $3$ & $2$ & $6$ & $9$  & $0.11$ & \num{-8.340e-01} & \num{ 1.776e-15}&32&70&0.08& $6$ & $39$  & $0.01$ & \num{-8.340e-01} & \num{ 5.446e-09} \\
HS40 & $4$ & $3$ & $4$ & $8$  & $0.06$ & \num{-2.500e-01} & \num{ 1.058e-12}&13&32&0.04& $3$ & $32$  & $0.03$ & \num{-2.500e-01} & \num{ 2.040e-09} \\
HS44 & $4$ & $6$ & $3$ & $8$  & $0.05$ & \num{-1.000e+00}(*) & \num{ 0.000e+00}&12&23&0.04& $4$ & $30$  & $0.02$ & \num{-1.500e+01} & \num{ 0.000e+00} \\
HS47 & $5$ & $3$ & $23$ & $28$  & $3.42$ & \num{ 1.989e-07} & \num{ 3.781e-11}&1391&7878&3.35& $72$ & $696$  & $0.38$ & \num{ 1.007e+01} & \num{ 6.661e-16} \\
HS50 & $5$ & $3$ & $19$ & $24$  & $1.11$ & \num{ 1.056e+00} & \num{ 2.398e-14}&126&423&1.06& $11$ & $85$  & $0.04$ & \num{ 1.339e-06} & \num{ 8.882e-16} \\
HS64 & $3$ & $1$ & $17$ & $20$  & $0.55$ & \num{ 6.300e+03} & \num{ 0.000e+00}&156&435&0.52& $7$ & $40$  & $0.02$ & \num{ 6.300e+03} & \num{ 7.051e-10} \\
HS66 & $3$ & $2$ & $2$ & $5$  & $0.03$ & \num{ 5.182e-01} & \num{ 3.563e-12}&7&12&0.02& $4$ & $27$  & $0.02$ & \num{ 5.182e-01} & \num{ 0.000e+00} \\
HS67 & $3$ & $14$ & $39$ & $42$  & $1.45$ & \num{-1.162e+03} & \num{ 0.000e+00}&450&1380&1.38& $10$ & $55$  & $0.03$ & \num{-1.162e+03} & \num{ 0.000e+00} \\
HS72 & $4$ & $2$ & $6$ & $10$  & $0.24$ & \num{ 7.277e+02} & \num{ 6.263e-14}&43&97&0.22& $4$ & $31$  & $0.02$ & \num{ 7.277e+02} & \num{ 1.572e-10} \\
HS75 & $4$ & $5$ & $16$ & $20$  & $0.31$ & \num{ 5.174e+03} & \num{ 5.684e-14}&97&180&0.27& $5$ & $37$  & $0.03$ & \num{ 5.174e+03} & \num{ 7.739e-08} \\
HS85 & $5$ & $21$ & $31$ & $36$  & $0.92$ & \num{-2.216e+00} & \num{ 0.000e+00}&357&903&0.83& $10$ & $78$  & $0.05$ & \num{-2.216e+00} & \num{ 2.389e-08} \\
HS87 & $6$ & $4$ & $2$ & $8$  & $0.11$ & \num{ 8.997e+03} & \num{ 7.844e-12}&15&43&0.1& $3$ & $33$  & $0.03$ & \num{ 8.997e+03} & \num{ 2.396e-08} \\
HS88 & $2$ & $1$ & $4$ & $6$  & $0.31$ & \num{ 1.363e+00} & \num{ 1.301e-14}&69&218&0.3& $23$ & $140$  & $0.08$ & \num{ 1.362e+00} & \num{ 4.623e-07} \\
HS89 & $3$ & $1$ & $5$ & $8$  & $0.55$ & \num{ 1.363e+00} & \num{ 0.000e+00}&113&432&0.54& $5$ & $35$  & $0.04$ & \num{ 7.283e-08}$(\dagger)$ & \num{ 1.332e-01} \\
HS90 & $4$ & $1$ & $11$ & $15$  & $3.6$ & \num{ 1.363e+00} & \num{ 8.118e-18}&718&3450&3.57& $24$ & $228$  & $0.1$ & \num{ 1.362e+00} & \num{ 5.472e-07} \\
HS93 & $6$ & $2$ & $25$ & $31$  & $1.57$ & \num{ 1.351e+02} & \num{ 3.331e-16}&371&1842&1.5& $8$ & $87$  & $0.05$ & \num{ 1.351e+02} & \num{ 0.000e+00} \\
HS98 & $6$ & $4$ & $2$ & $8$  & $0.03$ & \num{ 3.136e+00} & \num{ 0.000e+00}&4&6&0.02& $3$ & $32$  & $0.02$ & \num{ 3.136e+00} & \num{ 0.000e+00} \\
HS100 & $7$ & $4$ & $89$ & $96$  & $24.86$ & \num{ 6.806e+02} & \num{ 0.000e+00}&2860&19399&24.54& $14$ & $212$  & $0.07$ & \num{ 6.806e+02} & \num{ 0.000e+00} \\
HS101 & $7$ & $5$ & $71$ & $78$  & $12.74$ & \num{ 1.810e+03} & \num{ 2.776e-16}&2004&11919&12.48& $18$ & $198$  & $0.08$ & \num{ 1.810e+03} & \num{ 2.139e-08} \\
HS102 & $7$ & $5$ & $78$ & $85$  & $10.45$ & \num{ 9.119e+02} & \num{ 0.000e+00}&1601&8225&10.11& $15$ & $200$  & $0.07$ & \num{ 9.119e+02} & \num{ 5.354e-07} \\
HS103 & $7$ & $5$ & $39$ & $46$  & $4.22$ & \num{ 5.437e+02} & \num{ 6.106e-16}&685&3515&4.09& $16$ & $202$  & $0.09$ & \num{ 5.437e+02} & \num{ 1.040e-07} \\
HS104 & $8$ & $5$ & $80$ & $88$  & $13.31$ & \num{ 3.951e+00} & \num{ 0.000e+00}&2262&12381&13.0& $20$ & $259$  & $0.04$ & \num{ 3.951e+00} & \num{ 1.537e-07} \\
LOADBAL & $31$ & $31$ & $1$ & $15499$  & $0.88$ & \num{ 1.547e+00}(***) & \num{ 1.798e+308}&0&1&0.0& $17$ & $594$  & $0.09$ & \num{ 4.529e-01} & \num{ 1.483e-14} \\
OPTPRLOC & $30$ & $30$ & $8$ & $38$  & $0.52$ & \num{-1.642e+01} & \num{ 4.504e-10}&53&125&0.44& $13$ & $487$  & $0.11$ & \num{-1.642e+01} & \num{ 4.092e-08} \\
CB3 & $3$ & $3$ & $1$ & $4$  & $0.01$ & \num{ 2.000e+00} & \num{ 0.000e+00}&2&3&0.01& $1$ & $10$  & $0.0$ & \num{ 2.000e+00} & \num{ 7.270e-13} \\
CRESC50 & $6$ & $100$ & $1000$ & $1006$  & $214.64$ & \num{ 5.939e-01}(**) & \num{ 5.482e-13}&26160&96691&211.17& $438$ & $5787$  & $1.92$ & \num{ 5.946e-01} & \num{ 1.188e-14} \\
DEMBO7 & $16$ & $20$ & $19$ & $35$  & $2.06$ & \num{ 1.748e+02} & \num{ 4.073e-15}&302&795&1.97& $17$ & $354$  & $0.14$ & \num{ 1.749e+02} & \num{ 6.661e-16} \\
DNIEPER & $57$ & $24$ & $5$ & $62$  & $1.4$ & \num{ 1.874e+04} & \num{ 1.405e-09}&100&465&1.36& $2$ & $177$  & $0.06$ & \num{ 1.874e+04} & \num{ 1.800e-07} \\
EXPFITA & $5$ & $22$ & $1$ & $10999$  & $0.63$ & \num{ 2.999e+01}(***) & \num{ 1.798e+308}&0&1&0.0& $13$ & $98$  & $0.03$ & \num{ 1.137e-03} & \num{ 6.047e-12} \\
HIMMELBI & $100$ & $12$ & $114$ & $215$  & $68.51$ & \num{-1.467e+03}(*) & \num{ 1.789e-09}&2618&12586&67.29& $31$ & $3273$  & $0.41$ & \num{-1.736e+03} & \num{ 1.155e-14} \\
SYNTHES1 & $6$ & $6$ & $1$ & $2999$  & $0.21$ & \num{ 1.000e+01}(***) & \num{ 1.798e+308}&0&1&0.0& $7$ & $65$  & $0.02$ & \num{ 7.593e-01} & \num{ 3.261e-07} \\
TWOBARS & $2$ & $2$ & $3$ & $5$  & $0.06$ & \num{ 1.509e+00} & \num{ 4.742e-11}&16&55&0.05& $7$ & $49$  & $0.03$ & \num{ 1.509e+00} & \num{ 2.179e-09} \\
DIPIGRI & $7$ & $4$ & $87$ & $94$  & $78.71$ & \num{ 6.806e+02} & \num{ 1.265e-13}&7019&23030&78.4& $14$ & $212$  & $0.06$ & \num{ 6.806e+02} & \num{ 0.000e+00} \\
\end{tabular}}
\caption{Noiseless Problems with Feasible $x_0$. $\tau = 10^{-5}$; $n$: number of variables, $m$: number of constraints, \#iter: number of (outer) iterations, \#feval: number of function evaluations, time: total CPU time passed,  $f$: final objective value, feas err: final feasibility error, \#iter(sub): total number of iterations for solving TR subproblem, \#ceval: total number of constraint evaluations(note that the number of constraint evaluations and function evaluations are same for \texttt{KNITRO}), time(sub): total CPU time elapsed for solving subproblem. * indicates that \texttt{FIBO} terminates with singular interpolation system error, ** indicates that \texttt{FIBO} terminates with maximum number of iterations, *** indicates that \texttt{FIBO} terminates with maximum number of function evaluations. }
\label{tab:noiseless feas 1e-5}
\end{table}

%tau = 1e-7
\begin{table}[htp]
\centering
\resizebox{0.77\columnwidth}{!}{\begin{tabular}{lll | llllllll | lllll}
\hline
\multicolumn{3}{l}{} & \multicolumn{8}{c}{FIBO}       & \multicolumn{5}{c}{knitro}         \\ \hline
Problem   & n   & m  & \#iter & \#feval & time & $f$ & feas err & \#iter(sub) & \#ceval & time(sub) &  \#iter & \#feval & time & $f$ & feas err \\ \hline
HS13 & $2$ & $1$ & $3$ & $5$  & $0.29$ & \num{ 1.000e+00} & \num{ 1.309e-14}&136&191&0.28& $19$ & $100$  & $0.1$ & \num{ 9.983e-01} & \num{ 6.579e-10} \\
HS22 & $2$ & $2$ & $1$ & $3$  & $0.03$ & \num{ 1.000e+00} & \num{ 1.789e-12}&9&23&0.02& $0$ & $3$  & $0.0$ & \num{ 1.000e+00} & \num{ 1.843e-14} \\
HS23 & $2$ & $5$ & $3$ & $5$  & $0.04$ & \num{ 2.000e+00} & \num{ 0.000e+00}&11&14&0.03& $5$ & $24$  & $0.03$ & \num{ 2.000e+00} & \num{ 0.000e+00} \\
HS26 & $3$ & $1$ & $43$ & $46$  & $12.4$ & \num{ 2.878e-08} & \num{ 1.110e-16}&5194&38850&12.32& $28$ & $203$  & $0.08$ & \num{ 2.719e-11} & \num{ 8.518e-07} \\
HS32 & $3$ & $2$ & $4$ & $7$  & $0.05$ & \num{ 1.000e+00} & \num{ 0.000e+00}&10&16&0.03& $2$ & $15$  & $0.02$ & \num{ 1.000e+00} & \num{ 0.000e+00} \\
HS34 & $3$ & $2$ & $6$ & $9$  & $0.09$ & \num{-8.340e-01} & \num{ 1.776e-15}&32&70&0.07& $6$ & $39$  & $0.02$ & \num{-8.340e-01} & \num{ 5.446e-09} \\
HS40 & $4$ & $3$ & $4$ & $8$  & $0.05$ & \num{-2.500e-01} & \num{ 1.058e-12}&13&32&0.04& $3$ & $32$  & $0.05$ & \num{-2.500e-01} & \num{ 2.040e-09} \\
HS44 & $4$ & $6$ & $3$ & $8$  & $0.05$ & \num{-1.000e+00}(*) & \num{ 0.000e+00}&12&23&0.03& $4$ & $30$  & $0.02$ & \num{-1.500e+01} & \num{ 0.000e+00} \\
HS47 & $5$ & $3$ & $23$ & $28$  & $3.44$ & \num{ 1.989e-07} & \num{ 3.781e-11}&1391&7878&3.36& $72$ & $696$  & $0.36$ & \num{ 1.007e+01} & \num{ 6.661e-16} \\
HS50 & $5$ & $3$ & $19$ & $24$  & $1.12$ & \num{ 1.056e+00} & \num{ 2.398e-14}&126&423&1.06& $13$ & $100$  & $0.03$ & \num{ 2.358e-08} & \num{ 4.441e-16} \\
HS64 & $3$ & $1$ & $17$ & $20$  & $0.55$ & \num{ 6.300e+03} & \num{ 0.000e+00}&156&435&0.52& $7$ & $40$  & $0.02$ & \num{ 6.300e+03} & \num{ 7.051e-10} \\
HS66 & $3$ & $2$ & $2$ & $5$  & $0.03$ & \num{ 5.182e-01} & \num{ 3.563e-12}&7&12&0.02& $5$ & $33$  & $0.02$ & \num{ 5.182e-01} & \num{ 5.154e-11} \\
HS67 & $3$ & $14$ & $40$ & $43$  & $1.46$ & \num{-1.162e+03} & \num{ 2.274e-12}&458&1389&1.36& $11$ & $60$  & $0.03$ & \num{-1.162e+03} & \num{ 1.364e-12} \\
HS72 & $4$ & $2$ & $7$ & $11$  & $0.21$ & \num{ 7.277e+02} & \num{ 2.819e-18}&44&99&0.2& $4$ & $31$  & $0.02$ & \num{ 7.277e+02} & \num{ 1.572e-10} \\
HS75 & $4$ & $5$ & $16$ & $20$  & $0.31$ & \num{ 5.174e+03} & \num{ 5.684e-14}&97&180&0.28& $5$ & $37$  & $0.04$ & \num{ 5.174e+03} & \num{ 7.739e-08} \\
HS85 & $5$ & $21$ & $31$ & $36$  & $1.41$ & \num{-2.216e+00} & \num{ 0.000e+00}&357&903&1.3& $10$ & $78$  & $0.03$ & \num{-2.216e+00} & \num{ 2.389e-08} \\
HS87 & $6$ & $4$ & $9$ & $15$  & $0.7$ & \num{ 8.997e+03} & \num{ 2.274e-13}&80&287&0.67& $6$ & $73$  & $0.03$ & \num{ 8.997e+03} & \num{ 1.152e-08} \\
HS88 & $2$ & $1$ & $4$ & $6$  & $0.31$ & \num{ 1.363e+00} & \num{ 1.301e-14}&69&218&0.3& $23$ & $140$  & $0.07$ & \num{ 1.362e+00} & \num{ 4.623e-07} \\
HS89 & $3$ & $1$ & $6$ & $9$  & $0.68$ & \num{ 1.363e+00} & \num{ 5.535e-16}&142&532&0.66& $5$ & $35$  & $0.04$ & \num{ 7.283e-08}$(\dagger)$ & \num{ 1.332e-01} \\
HS90 & $4$ & $1$ & $11$ & $15$  & $3.6$ & \num{ 1.363e+00} & \num{ 8.118e-18}&718&3450&3.57& $24$ & $228$  & $0.1$ & \num{ 1.362e+00} & \num{ 5.472e-07} \\
HS93 & $6$ & $2$ & $36$ & $42$  & $2.18$ & \num{ 1.351e+02} & \num{ 0.000e+00}&483&2291&2.06& $9$ & $96$  & $0.05$ & \num{ 1.351e+02} & \num{ 0.000e+00} \\
HS98 & $6$ & $4$ & $2$ & $8$  & $0.02$ & \num{ 3.136e+00} & \num{ 0.000e+00}&4&6&0.02& $3$ & $32$  & $0.01$ & \num{ 3.136e+00} & \num{ 0.000e+00} \\
HS100 & $7$ & $4$ & $106$ & $113$  & $28.46$ & \num{ 6.806e+02} & \num{ 0.000e+00}&3280&22448&28.08& $28$ & $362$  & $0.19$ & \num{ 6.806e+02} & \num{ 5.093e-11} \\
HS101 & $7$ & $5$ & $90$ & $97$  & $15.48$ & \num{ 1.810e+03} & \num{ 2.742e-12}&2377&14118&15.15& $20$ & $216$  & $0.1$ & \num{ 1.810e+03} & \num{ 1.032e-08} \\
HS102 & $7$ & $5$ & $89$ & $96$  & $12.09$ & \num{ 9.119e+02} & \num{ 2.012e-16}&1847&9679&11.74& $15$ & $200$  & $0.08$ & \num{ 9.119e+02} & \num{ 5.354e-07} \\
HS103 & $7$ & $5$ & $42$ & $49$  & $4.08$ & \num{ 5.437e+02} & \num{ 4.302e-16}&713&3618&3.92& $16$ & $202$  & $0.08$ & \num{ 5.437e+02} & \num{ 1.040e-07} \\
HS104 & $8$ & $5$ & $88$ & $96$  & $13.71$ & \num{ 3.951e+00} & \num{ 4.441e-16}&2400&13125&13.38& $20$ & $259$  & $0.05$ & \num{ 3.951e+00} & \num{ 1.537e-07} \\
LOADBAL & $31$ & $31$ & $1$ & $15499$  & $0.89$ & \num{ 1.547e+00}(***) & \num{ 1.798e+308}&0&1&0.0& $18$ & $627$  & $0.09$ & \num{ 4.529e-01} & \num{ 8.549e-15} \\
OPTPRLOC & $30$ & $30$ & $8$ & $38$  & $0.51$ & \num{-1.642e+01} & \num{ 4.504e-10}&53&125&0.43& $13$ & $487$  & $0.11$ & \num{-1.642e+01} & \num{ 4.092e-08} \\
CB3 & $3$ & $3$ & $1$ & $4$  & $0.01$ & \num{ 2.000e+00} & \num{ 0.000e+00}&2&3&0.0& $1$ & $10$  & $0.0$ & \num{ 2.000e+00} & \num{ 7.270e-13} \\
CRESC50 & $6$ & $100$ & $1000$ & $1006$  & $214.17$ & \num{ 5.939e-01}(**) & \num{ 5.482e-13}&26160&96691&210.8& $438$ & $5787$  & $1.87$ & \num{ 5.946e-01} & \num{ 1.188e-14} \\
DEMBO7 & $16$ & $20$ & $19$ & $35$  & $2.11$ & \num{ 1.748e+02} & \num{ 4.073e-15}&302&795&2.04& $17$ & $354$  & $0.15$ & \num{ 1.749e+02} & \num{ 6.661e-16} \\
DNIEPER & $57$ & $24$ & $23$ & $80$  & $3.43$ & \num{ 1.874e+04} & \num{ 5.684e-14}&272&1088&3.27& $6$ & $416$  & $0.08$ & \num{ 1.874e+04} & \num{ 5.623e-07} \\
EXPFITA & $5$ & $22$ & $1$ & $10999$  & $0.6$ & \num{ 2.999e+01}(***) & \num{ 1.798e+308}&0&1&0.0& $14$ & $105$  & $0.04$ & \num{ 1.137e-03} & \num{ 1.776e-14} \\
HIMMELBI & $100$ & $12$ & $114$ & $215$  & $69.11$ & \num{-1.467e+03}(*) & \num{ 1.789e-09}&2618&12586&68.01& $41$ & $4300$  & $0.52$ & \num{-1.736e+03} & \num{ 7.105e-15} \\
SYNTHES1 & $6$ & $6$ & $1$ & $2999$  & $0.23$ & \num{ 1.000e+01}(***) & \num{ 1.798e+308}&0&1&0.0& $7$ & $65$  & $0.02$ & \num{ 7.593e-01} & \num{ 3.261e-07} \\
TWOBARS & $2$ & $2$ & $5$ & $7$  & $0.1$ & \num{ 1.509e+00} & \num{ 0.000e+00}&21&84&0.08& $7$ & $49$  & $0.03$ & \num{ 1.509e+00} & \num{ 2.179e-09} \\
DIPIGRI & $7$ & $4$ & $96$ & $103$  & $80.88$ & \num{ 6.806e+02} & \num{ 4.547e-13}&7235&24605&80.56& $28$ & $365$  & $0.1$ & \num{ 6.806e+02} & \num{ 7.237e-11} \\
\end{tabular}}
\caption{Noiseless Problems with Feasible $x_0$. $\tau = 10^{-7}$; $n$: number of variables, $m$: number of constraints, \#iter: number of (outer) iterations, \#feval: number of function evaluations, time: total CPU time passed,  $f$: final objective value, feas err: final feasibility error, \#iter(sub): total number of iterations for solving TR subproblem, \#ceval: total number of constraint evaluations(note that the number of constraint evaluations and function evaluations are same for \texttt{KNITRO}), time(sub): total CPU time elapsed for solving subproblem. * indicates that \texttt{FIBO} terminates with singular interpolation system error, ** indicates that \texttt{FIBO} terminates with maximum number of iterations, *** indicates that \texttt{FIBO} terminates with maximum number of function evaluations. }
\label{tab:noiseless feas 1e-7}
\end{table}

\end{landscape}

\section{Numerical Results for Infeasible Initial Point}
\label{app: outerdfo infeas}
We now include the cost of obtaining an initial feasible point for \texttt{FIBO}. In other words, the CPU time and number of constraint evaluations associated with obtaining a feasible starting point are counted. Simple comparison with the previous table indicates that only a few constraint evaluations are needed in order to obtain a feasible solution. As for \texttt{KNITRO}, we instead run from the given initial point $x_0$, which is potentially infeasible.

\begin{figure}[htp]
	\centering
	\includegraphics[width = 0.48\textwidth]{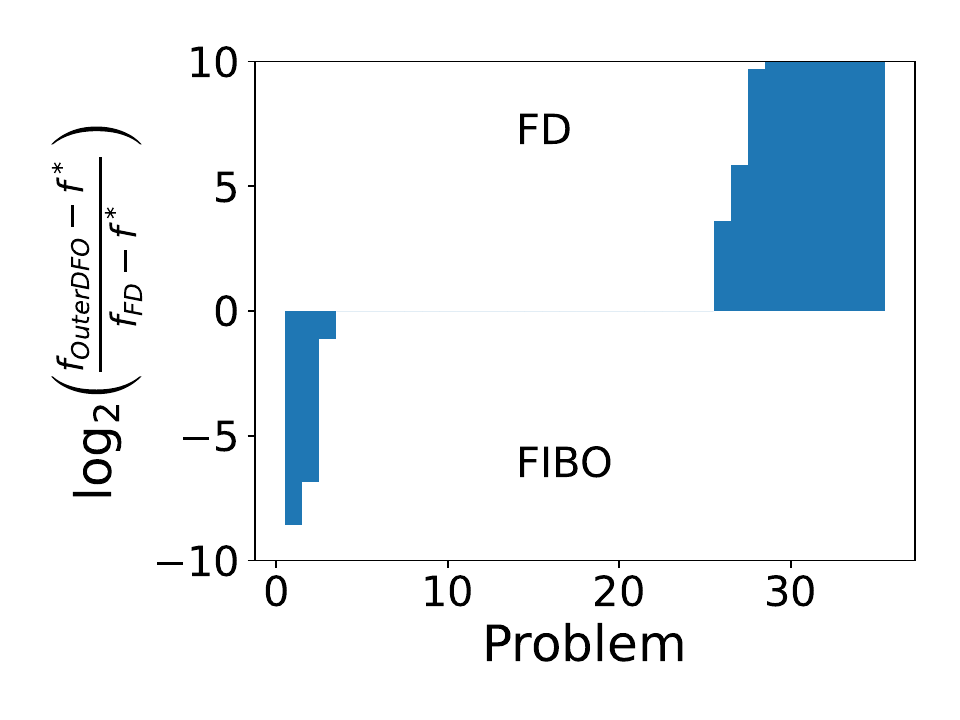}
	\caption{Infeasible $x_0$. Log-ratio Plot for Comparing the Final Accuracy \eqref{acc} of \texttt{FIBO} and \texttt{FD}.}
	\label{fig:infeas acc}
\end{figure}

\begin{figure}[htp]
	\centering
	\includegraphics[width=0.38\textwidth]{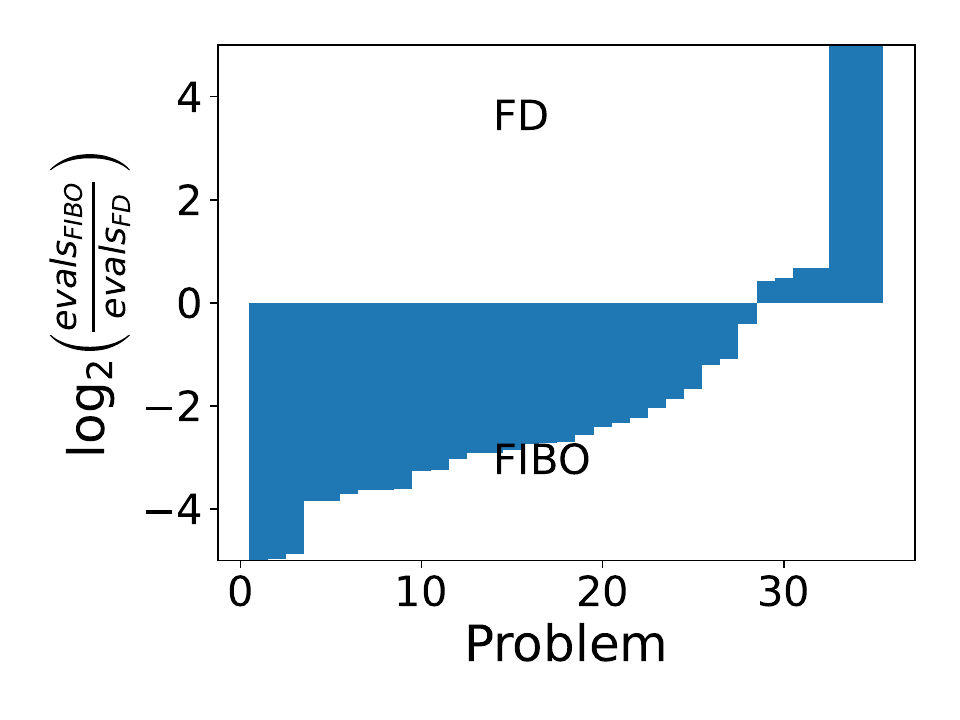}
	\includegraphics[width=0.38\textwidth]{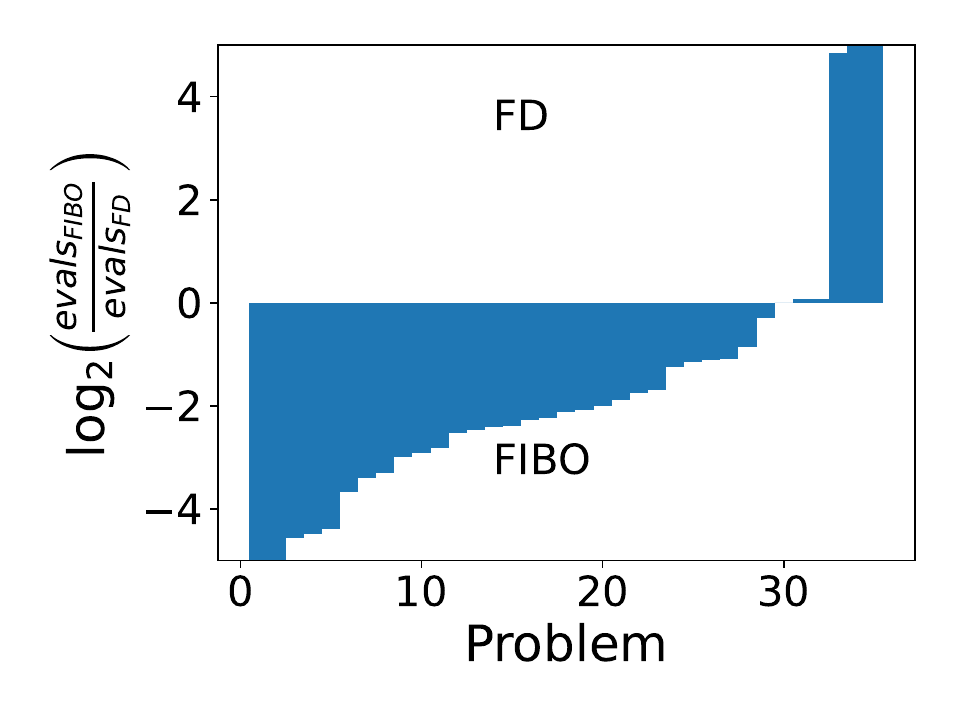} \\
	\includegraphics[width=0.38\textwidth]{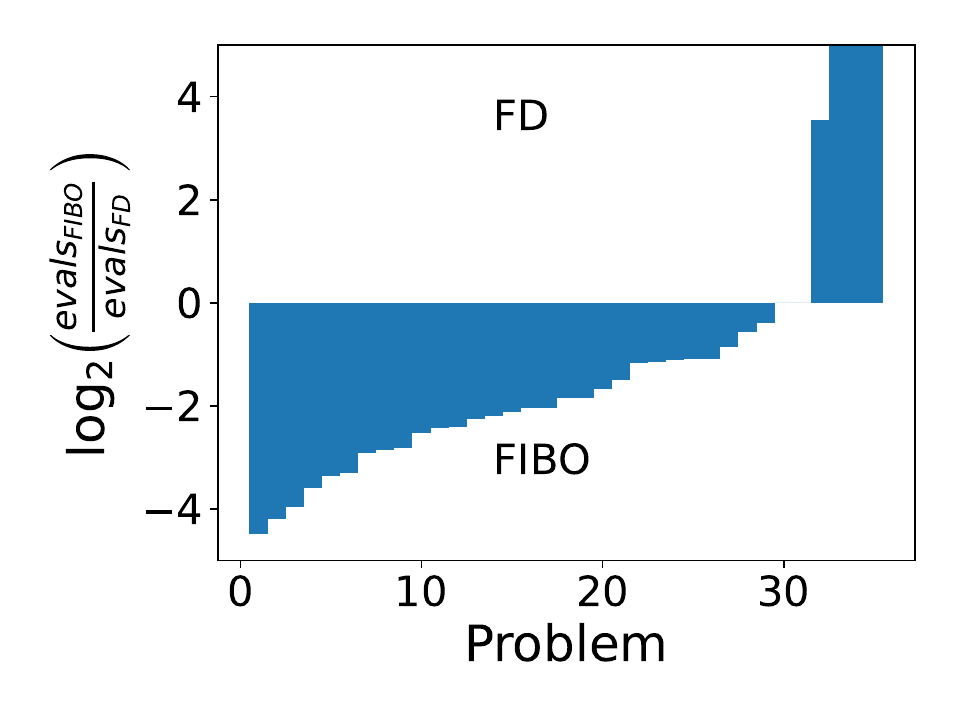}
	\includegraphics[width=0.38\textwidth]{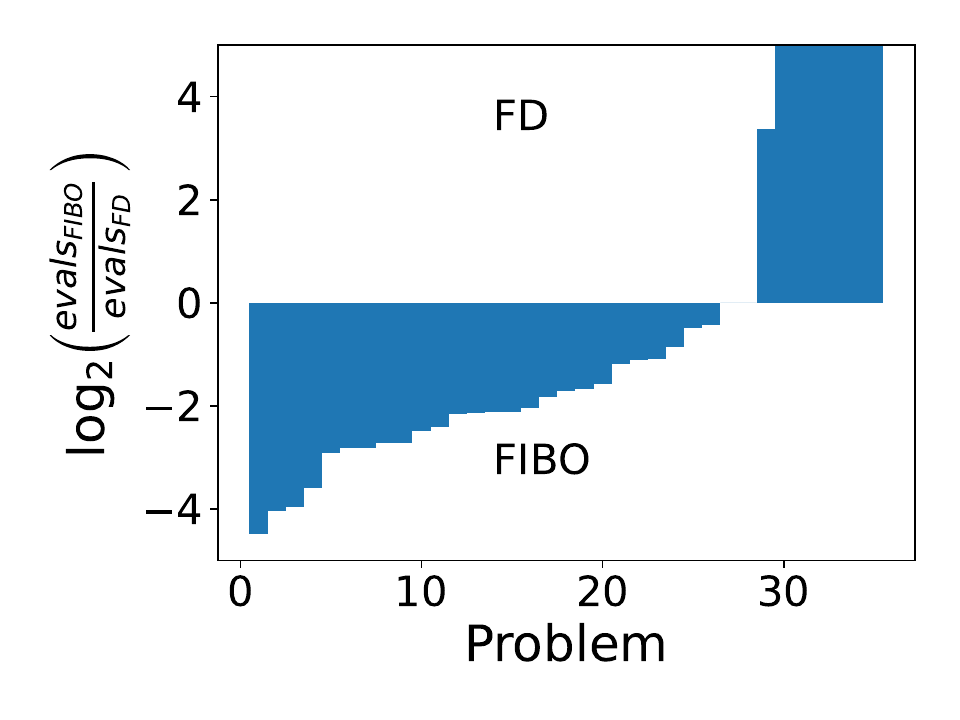}
	\caption{Infeasible $x_0$. Log-ratio plot comparing \texttt{FIBO} and \texttt{FD} in terms of the number of function evaluations to satisfy \eqref{efficiency} for $\tau = 10^{-1} \text{ (upper left)}, 10^{-3}\text{ (upper right)}, 10^{-5}\text{ (bottom left)}, 10^{-7}\text{ (bottom right)}$.}
	\label{fig:infeas feval}
\end{figure}

\begin{figure}[htp]
	\centering
	\includegraphics[width=0.38\textwidth]{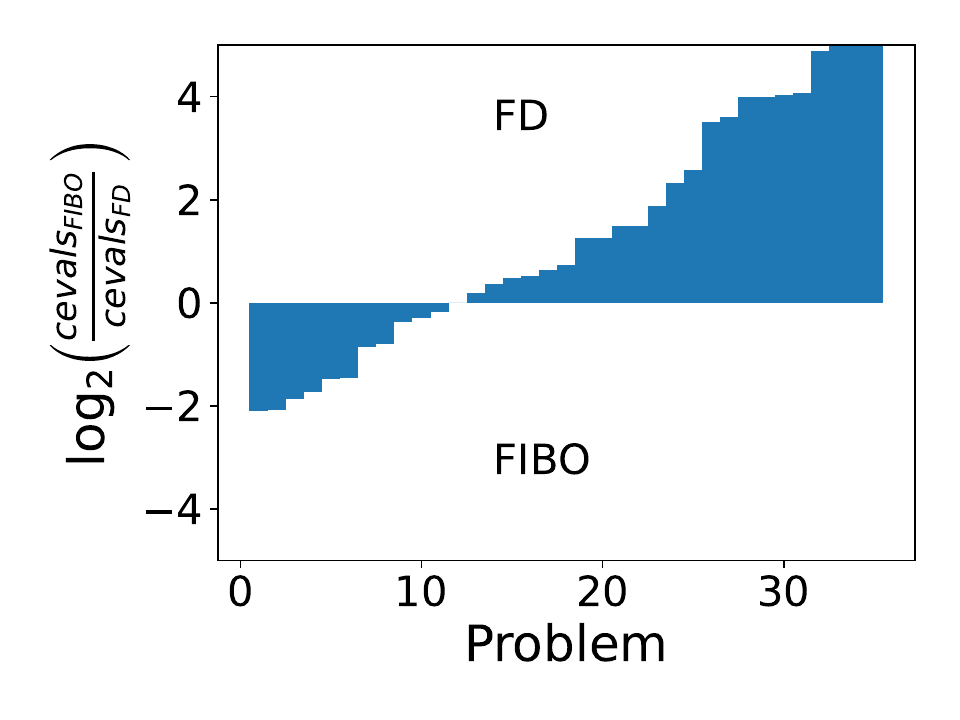}
	\includegraphics[width=0.38\textwidth]{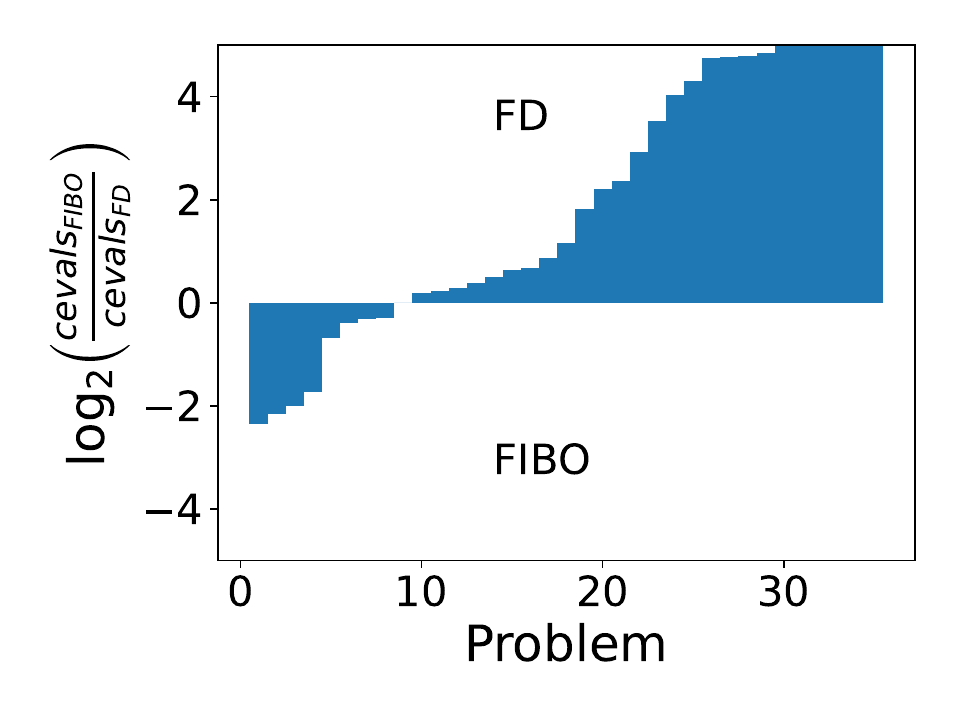} \\
	\includegraphics[width=0.38\textwidth]{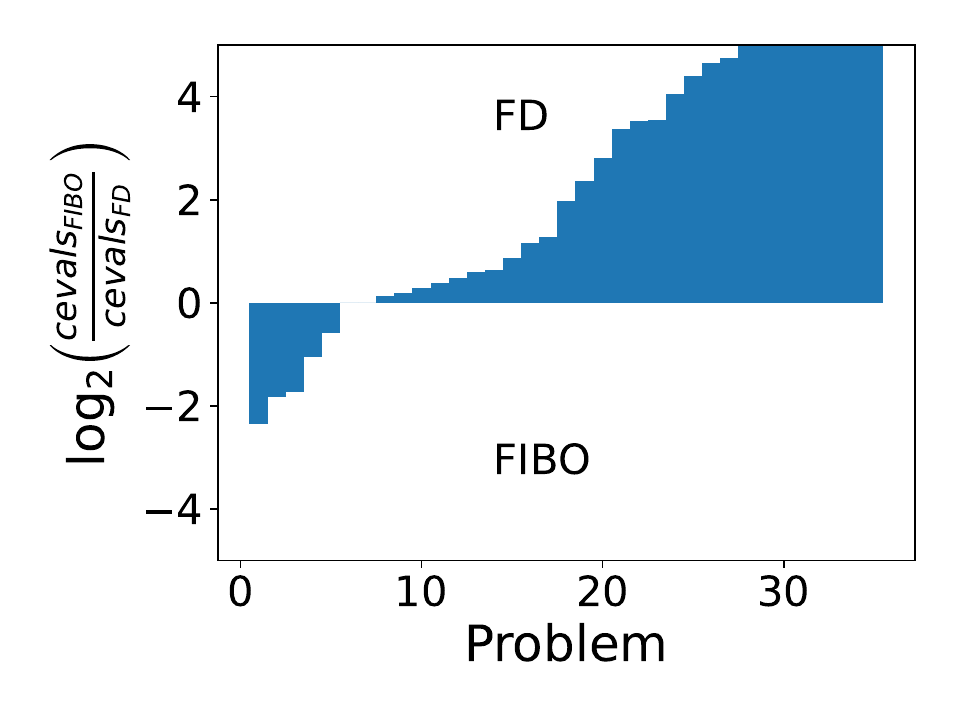}
	\includegraphics[width=0.38\textwidth]{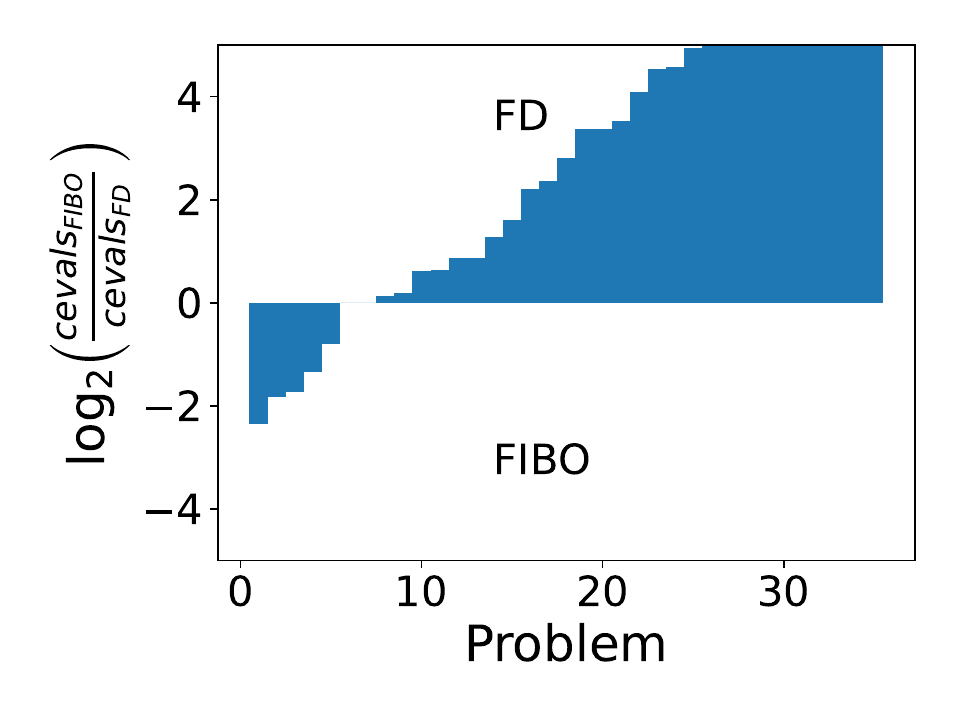}
	\caption{Infeasible $x_0$. Log-ratio plot comparing \texttt{FIBO} and \texttt{FD} in terms of the number of constraint evaluations to satisfy \eqref{efficiency} for $\tau = 10^{-1} \text{ (upper left)}, 10^{-3}\text{ (upper right)}, 10^{-5}\text{ (bottom left)}, 10^{-7}\text{ (bottom right)}$.}
	\label{fig:infeas ceval}
\end{figure}

As indicated by Figures \ref{fig:infeas acc}-\ref{fig:infeas ceval}, the conclusions remain the same as in Section \ref{sec:num} and the cost of obtaining a feasible starting point is quite negligible.

\begin{landscape}
\begin{table}[]
\centering
\resizebox{0.77\columnwidth}{!}{\begin{tabular}{lll | llllllll | lllll}
\hline
\multicolumn{3}{l}{} & \multicolumn{8}{c}{FIBO}       & \multicolumn{5}{c}{knitro}         \\ \hline
Problem   & n   & m  & \#iter & \#feval & time & $f$ & feas err & \#iter(sub) & \#ceval & time(sub) &  \#iter & \#feval & time & $f$ & feas err \\ \hline
HS13 & $2$ & $1$ & $3$ & $5$  & $0.32$ & \num{ 1.000e+00} & \num{ 1.309e-14}&136&193&0.28& $20$ & $106$  & $0.09$ & \num{ 9.983e-01} & \num{ 6.579e-10} \\
HS22 & $2$ & $2$ & $2$ & $4$  & $0.05$ & \num{ 1.000e+00} & \num{ 0.000e+00}&10&27&0.03& $5$ & $24$  & $0.01$ & \num{ 1.000e+00} & \num{ 0.000e+00} \\
HS23 & $2$ & $5$ & $4$ & $6$  & $0.05$ & \num{ 2.000e+00} & \num{ 0.000e+00}&12&18&0.03& $7$ & $32$  & $0.01$ & \num{ 2.000e+00} & \num{ 0.000e+00} \\
HS26 & $3$ & $1$ & $53$ & $56$  & $13.0$ & \num{ 7.103e-10} & \num{ 2.220e-16}&5245&39138&12.85& $34$ & $244$  & $0.11$ & \num{ 7.434e-12} & \num{ 6.661e-16} \\
HS32 & $3$ & $2$ & $4$ & $7$  & $0.05$ & \num{ 1.000e+00} & \num{ 0.000e+00}&10&17&0.03& $2$ & $15$  & $0.01$ & \num{ 1.000e+00} & \num{ 0.000e+00} \\
HS34 & $3$ & $2$ & $6$ & $9$  & $0.12$ & \num{-8.340e-01} & \num{ 1.776e-15}&32&71&0.08& $7$ & $44$  & $0.02$ & \num{-8.340e-01} & \num{ 0.000e+00} \\
HS40 & $4$ & $3$ & $22$ & $26$  & $0.37$ & \num{-2.500e-01} & \num{ 1.110e-16}&56&229&0.26& $7$ & $51$  & $0.02$ & \num{-2.500e-01} & \num{ 1.665e-16} \\
HS44 & $4$ & $6$ & $3$ & $8$  & $0.06$ & \num{-1.000e+00}(*) & \num{ 0.000e+00}&12&24&0.03& $5$ & $36$  & $0.01$ & \num{-1.500e+01} & \num{ 0.000e+00} \\
HS47 & $5$ & $3$ & $23$ & $28$  & $3.67$ & \num{ 1.989e-07} & \num{ 3.781e-11}&1391&7879&3.47& $72$ & $696$  & $0.35$ & \num{ 1.007e+01} & \num{ 6.661e-16} \\
HS50 & $5$ & $3$ & $19$ & $24$  & $1.18$ & \num{ 1.056e+00} & \num{ 2.398e-14}&126&424&1.1& $16$ & $121$  & $0.04$ & \num{ 7.669e-17} & \num{ 4.441e-16} \\
HS64 & $3$ & $1$ & $43$ & $46$  & $1.07$ & \num{ 6.300e+03} & \num{ 0.000e+00}&300&674&0.95& $16$ & $89$  & $0.04$ & \num{ 6.300e+03} & \num{ 1.110e-16} \\
HS66 & $3$ & $2$ & $2$ & $6$  & $0.04$ & \num{ 5.182e-01}(*) & \num{ 3.563e-12}&7&13&0.02& $7$ & $43$  & $0.01$ & \num{ 5.182e-01} & \num{ 0.000e+00} \\
HS67 & $3$ & $14$ & $45$ & $49$  & $1.71$ & \num{-1.162e+03}(*) & \num{ 0.000e+00}&472&1409&1.59& $14$ & $75$  & $0.03$ & \num{-1.162e+03} & \num{ 0.000e+00} \\
HS72 & $4$ & $2$ & $7$ & $11$  & $0.32$ & \num{ 7.277e+02} & \num{ 2.819e-18}&44&111&0.21& $15$ & $96$  & $0.05$ & \num{ 7.277e+02} & \num{ 1.301e-18} \\
HS75 & $4$ & $5$ & $16$ & $20$  & $0.54$ & \num{ 5.174e+03} & \num{ 5.684e-14}&97&185&0.43& $7$ & $48$  & $0.03$ & \num{ 5.174e+03} & \num{ 2.842e-14} \\
HS85 & $5$ & $21$ & $31$ & $36$  & $1.47$ & \num{-2.216e+00} & \num{ 0.000e+00}&357&904&1.31& $11$ & $85$  & $0.03$ & \num{-2.216e+00} & \num{ 1.421e-14} \\
HS87 & $6$ & $4$ & $9$ & $15$  & $0.76$ & \num{ 8.997e+03} & \num{ 2.274e-13}&80&290&0.71& $14$ & $216$  & $0.07$ & \num{ 8.997e+03} & \num{ 2.274e-13} \\
HS88 & $2$ & $1$ & $5$ & $7$  & $0.43$ & \num{ 1.363e+00} & \num{ 0.000e+00}&92&329&0.39& $19$ & $106$  & $0.06$ & \num{ 1.363e+00} & \num{ 0.000e+00} \\
HS89 & $3$ & $1$ & $10$ & $13$  & $1.52$ & \num{ 1.363e+00} & \num{ 0.000e+00}&246&1033&1.39& $28$ & $172$  & $0.09$ & \num{ 1.363e+00} & \num{ 9.636e-18} \\
HS90 & $4$ & $1$ & $12$ & $16$  & $4.26$ & \num{ 1.363e+00} & \num{ 0.000e+00}&738&3541&4.18& $48$ & $426$  & $0.24$ & \num{ 1.363e+00} & \num{ 0.000e+00} \\
HS93 & $6$ & $2$ & $42$ & $48$  & $2.9$ & \num{ 1.351e+02} & \num{ 4.800e-12}&517&2434&2.61& $26$ & $295$  & $0.08$ & \num{ 1.351e+02} & \num{ 0.000e+00} \\
HS98 & $6$ & $4$ & $3$ & $10$  & $0.09$ & \num{ 3.136e+00}(*) & \num{ 1.776e-15}&5&13&0.02& $6$ & $56$  & $0.01$ & \num{ 3.136e+00} & \num{ 0.000e+00} \\
HS100 & $7$ & $4$ & $106$ & $113$  & $29.47$ & \num{ 6.806e+02} & \num{ 0.000e+00}&3280&22449&28.72& $51$ & $587$  & $0.16$ & \num{ 6.806e+02} & \num{ 7.161e-15} \\
HS101 & $7$ & $5$ & $103$ & $110$  & $19.45$ & \num{ 1.810e+03} & \num{ 3.955e-16}&2616&15705&18.39& $56$ & $582$  & $0.22$ & \num{ 1.810e+03} & \num{ 0.000e+00} \\
HS102 & $7$ & $5$ & $103$ & $110$  & $15.14$ & \num{ 9.119e+02} & \num{ 3.816e-16}&1985&10248&13.86& $37$ & $386$  & $0.13$ & \num{ 9.119e+02} & \num{ 2.498e-16} \\
HS103 & $7$ & $5$ & $56$ & $63$  & $5.53$ & \num{ 5.437e+02} & \num{ 1.318e-16}&825&4079&4.59& $28$ & $284$  & $0.08$ & \num{ 5.437e+02} & \num{ 2.220e-16} \\
HS104 & $8$ & $5$ & $245$ & $253$  & $43.52$ & \num{ 3.951e+00} & \num{ 9.714e-17}&6813&37274&41.25& $17$ & $190$  & $0.04$ & \num{ 3.951e+00} & \num{ 1.943e-16} \\
LOADBAL & $31$ & $31$ & $1$ & $15499$  & $0.94$ & \num{ 1.547e+00}(***) & \num{ 1.798e+308}&0&2&0.0& $40$ & $1380$  & $0.19$ & \num{ 4.529e-01} & \num{ 2.665e-15} \\
OPTPRLOC & $30$ & $30$ & $9$ & $39$  & $1.01$ & \num{-1.642e+01} & \num{ 1.990e-13}&54&130&0.45& $12$ & $456$  & $0.06$ & \num{-1.642e+01} & \num{ 1.279e-13} \\
CB3 & $3$ & $3$ & $1$ & $4$  & $0.03$ & \num{ 2.000e+00} & \num{ 0.000e+00}&2&9&0.0& $7$ & $40$  & $0.01$ & \num{ 2.000e+00} & \num{ 0.000e+00} \\
CRESC50 & $6$ & $100$ & $1000$ & $1006$  & $215.05$ & \num{ 5.939e-01}(**) & \num{ 5.482e-13}&26160&96697&210.07& $456$ & $6023$  & $2.46$ & \num{ 5.948e-01} & \num{ 4.547e-13} \\
DEMBO7 & $16$ & $20$ & $19$ & $35$  & $2.42$ & \num{ 1.748e+02} & \num{ 4.073e-15}&302&800&2.07& $45$ & $893$  & $0.2$ & \num{ 1.748e+02} & \num{ 3.608e-16} \\
DNIEPER & $57$ & $24$ & $24$ & $81$  & $7.75$ & \num{ 1.874e+04} & \num{ 5.684e-14}&273&1097&3.25& $6$ & $418$  & $0.07$ & \num{ 1.874e+04} & \num{ 5.684e-14} \\
EXPFITA & $5$ & $22$ & $1$ & $10999$  & $0.62$ & \num{ 2.999e+01}(***) & \num{ 1.798e+308}&0&2&0.0& $19$ & $140$  & $0.05$ & \num{ 1.137e-03} & \num{ 1.421e-14} \\
HIMMELBI & $100$ & $12$ & $114$ & $215$  & $72.25$ & \num{-1.467e+03}(*) & \num{ 1.789e-09}&2618&12588&67.95& $96$ & $9950$  & $1.01$ & \num{-1.736e+03} & \num{ 2.842e-14} \\
SYNTHES1 & $6$ & $6$ & $1$ & $2999$  & $0.24$ & \num{ 1.000e+01}(***) & \num{ 1.798e+308}&0&2&0.0& $9$ & $81$  & $0.02$ & \num{ 7.593e-01} & \num{ 1.110e-16} \\
TWOBARS & $2$ & $2$ & $7$ & $9$  & $0.13$ & \num{ 1.509e+00} & \num{ 0.000e+00}&24&113&0.09& $11$ & $57$  & $0.02$ & \num{ 1.509e+00} & \num{ 2.220e-16} \\
DIPIGRI & $7$ & $4$ & $116$ & $123$  & $85.87$ & \num{ 6.806e+02} & \num{ 3.581e-11}&7461&26265&84.59& $48$ & $581$  & $0.16$ & \num{ 6.806e+02} & \num{ 0.000e+00} \\
\end{tabular}}
\caption{Noiseless Problems with Infeasible $x_0$; $n$: number of variables, $m$: number of constraints, \#iter: number of (outer) iterations, \#feval: number of function evaluations, time: total CPU time passed,  $f$: final objective value, feas err: final feasibility error, \#iter(sub): total number of iterations for solving TR subproblem, \#ceval: total number of constraint evaluations(note that the number of constraint evaluations and function evaluations are same for \texttt{KNITRO}), time(sub): total CPU time elapsed for solving subproblem. * indicates that \texttt{FIBO} terminates with singular interpolation system error, ** indicates that \texttt{FIBO} terminates with maximum number of iterations, *** indicates that \texttt{FIBO} terminates with maximum number of function evaluations. }
\label{tab: noiseless acc infeas}
\end{table}

\end{landscape}

%tau = 1e-1
\begin{landscape}
\begin{table}[htp]
\centering
\resizebox{0.77\columnwidth}{!}{\begin{tabular}{lll | llllllll | lllll}
\hline
\multicolumn{3}{l}{} & \multicolumn{8}{c}{FIBO}       & \multicolumn{5}{c}{knitro}         \\ \hline
Problem   & n   & m  & \#iter & \#feval & time & $f$ & feas err & \#iter(sub) & \#ceval & time(sub) &  \#iter & \#feval & time & $f$ & feas err \\ \hline
HS13 & $2$ & $1$ & $2$ & $4$  & $0.2$ & \num{ 1.000e+00} & \num{ 1.310e-14}&75&106&0.19& $8$ & $38$  & $0.02$ & \num{ 1.080e+00} & \num{ 0.000e+00} \\
HS22 & $2$ & $2$ & $1$ & $3$  & $0.04$ & \num{ 1.000e+00} & \num{ 1.789e-12}&9&25&0.02& $3$ & $16$  & $0.01$ & \num{ 1.000e+00} & \num{ 1.231e-06} \\
HS23 & $2$ & $5$ & $3$ & $5$  & $0.07$ & \num{ 2.000e+00} & \num{ 0.000e+00}&11&16&0.03& $3$ & $16$  & $0.01$ & \num{ 2.075e+00} & \num{ 0.000e+00} \\
HS26 & $3$ & $1$ & $4$ & $7$  & $0.09$ & \num{ 1.996e-02} & \num{ 9.592e-14}&32&67&0.07& $14$ & $87$  & $0.03$ & \num{ 4.252e-07} & \num{ 9.633e-07} \\
HS32 & $3$ & $2$ & $4$ & $7$  & $0.1$ & \num{ 1.000e+00} & \num{ 0.000e+00}&10&17&0.04& $2$ & $15$  & $0.01$ & \num{ 1.000e+00} & \num{ 0.000e+00} \\
HS34 & $3$ & $2$ & $5$ & $8$  & $0.11$ & \num{-7.820e-01} & \num{ 1.066e-14}&31&69&0.08& $4$ & $29$  & $0.01$ & \num{-8.266e-01} & \num{ 0.000e+00} \\
HS40 & $4$ & $3$ & $1$ & $5$  & $0.03$ & \num{-2.454e-01} & \num{ 2.554e-15}&5&19&0.01& $4$ & $33$  & $0.01$ & \num{-2.500e-01} & \num{ 5.162e-08} \\
HS44 & $4$ & $6$ & $3$ & $8$  & $0.1$ & \num{-1.000e+00}(*) & \num{ 0.000e+00}&12&24&0.07& $4$ & $30$  & $0.01$ & \num{-1.500e+01} & \num{ 0.000e+00} \\
HS47 & $5$ & $3$ & $6$ & $11$  & $0.15$ & \num{ 5.639e-02} & \num{ 1.466e-12}&46&119&0.12& $17$ & $134$  & $0.04$ & \num{ 4.638e-08} & \num{ 3.167e-09} \\
HS50 & $5$ & $3$ & $19$ & $24$  & $1.26$ & \num{ 1.056e+00} & \num{ 2.398e-14}&126&423&1.19& $9$ & $71$  & $0.03$ & \num{ 4.941e-02} & \num{ 4.441e-16} \\
HS64 & $3$ & $1$ & $1$ & $4$  & $0.09$ & \num{ 6.813e+03} & \num{ 4.742e-11}&14&84&0.05& $5$ & $30$  & $0.02$ & \num{ 6.705e+03} & \num{ 3.942e-05} \\
HS66 & $3$ & $2$ & $1$ & $4$  & $0.05$ & \num{ 5.800e-01} & \num{ 3.563e-12}&6&11&0.01& $0$ & $3$  & $0.0$ & \num{ 5.800e-01} & \num{ 0.000e+00} \\
HS67 & $3$ & $14$ & $31$ & $34$  & $1.92$ & \num{-1.053e+03} & \num{ 0.000e+00}&402&1320&1.82& $8$ & $45$  & $0.02$ & \num{-1.157e+03} & \num{ 0.000e+00} \\
HS72 & $4$ & $2$ & $1$ & $5$  & $0.03$ & \num{ 7.407e+02} & \num{ 0.000e+00}&2&17&0.01& $11$ & $72$  & $0.04$ & \num{ 7.277e+02} & \num{ 8.766e-07} \\
HS75 & $4$ & $5$ & $1$ & $5$  & $0.06$ & \num{ 5.350e+03} & \num{ 1.259e-08}&4&13&0.01& $5$ & $36$  & $0.03$ & \num{ 5.174e+03} & \num{ 1.880e-04} \\
HS85 & $5$ & $21$ & $29$ & $34$  & $1.23$ & \num{-2.061e+00} & \num{ 1.421e-14}&348&881&1.04& $10$ & $78$  & $0.03$ & \num{-2.216e+00} & \num{ 2.389e-08} \\
HS87 & $6$ & $4$ & $1$ & $7$  & $0.1$ & \num{ 8.997e+03} & \num{ 2.274e-13}&9&27&0.06& $3$ & $33$  & $0.02$ & \num{ 8.997e+03} & \num{ 3.775e-08} \\
HS88 & $2$ & $1$ & $1$ & $3$  & $0.13$ & \num{ 1.385e+00} & \num{ 2.158e-14}&19&52&0.07& $16$ & $94$  & $0.05$ & \num{ 1.363e+00} & \num{ 1.074e-08} \\
HS89 & $3$ & $1$ & $2$ & $5$  & $0.41$ & \num{ 1.365e+00} & \num{ 1.123e-16}&57&349&0.28& $23$ & $147$  & $0.07$ & \num{ 1.362e+00} & \num{ 7.400e-07} \\
HS90 & $4$ & $1$ & $2$ & $6$  & $0.56$ & \num{ 1.446e+00} & \num{ 6.661e-16}&86&464&0.52& $36$ & $334$  & $0.15$ & \num{ 1.363e+00} & \num{ 1.142e-07} \\
HS93 & $6$ & $2$ & $1$ & $7$  & $0.13$ & \num{ 1.371e+02} & \num{ 2.946e-11}&10&25&0.11& $0$ & $5$  & $0.0$ & \num{ 1.371e+02} & \num{ 0.000e+00} \\
HS98 & $6$ & $4$ & $2$ & $8$  & $0.06$ & \num{ 3.136e+00} & \num{ 0.000e+00}&4&11&0.02& $4$ & $40$  & $0.01$ & \num{ 3.214e+00} & \num{ 0.000e+00} \\
HS100 & $7$ & $4$ & $1$ & $8$  & $0.14$ & \num{ 7.140e+02} & \num{ 0.000e+00}&22&80&0.12& $0$ & $5$  & $0.0$ & \num{ 7.140e+02} & \num{ 0.000e+00} \\
HS101 & $7$ & $5$ & $26$ & $33$  & $5.28$ & \num{ 1.982e+03} & \num{ 3.053e-16}&740&4391&5.04& $26$ & $269$  & $0.1$ & \num{ 1.810e+03} & \num{ 2.249e-05} \\
HS102 & $7$ & $5$ & $34$ & $41$  & $4.59$ & \num{ 9.948e+02} & \num{ 2.776e-17}&641&3328&4.29& $25$ & $272$  & $0.09$ & \num{ 9.129e+02} & \num{ 5.866e-05} \\
HS103 & $7$ & $5$ & $15$ & $22$  & $1.72$ & \num{ 5.623e+02} & \num{ 8.913e-12}&276&1270&1.63& $20$ & $212$  & $0.07$ & \num{ 5.434e+02} & \num{ 3.104e-04} \\
HS104 & $8$ & $5$ & $1$ & $9$  & $0.15$ & \num{ 4.200e+00} & \num{ 1.459e-11}&14&47&0.11& $11$ & $130$  & $0.03$ & \num{ 3.951e+00} & \num{ 5.769e-07} \\
LOADBAL & $31$ & $31$ & $1$ & $15499$  & $0.99$ & \num{ 1.547e+00}(***) & \num{ 1.798e+308}&0&2&0.0& $6$ & $231$  & $0.04$ & \num{ 5.252e-01} & \num{ 1.710e-14} \\
OPTPRLOC & $30$ & $30$ & $5$ & $35$  & $0.48$ & \num{-1.567e+01} & \num{ 2.220e-16}&44&107&0.4& $12$ & $456$  & $0.06$ & \num{-1.642e+01} & \num{ 1.279e-13} \\
CB3 & $3$ & $3$ & $1$ & $4$  & $0.02$ & \num{ 2.000e+00} & \num{ 0.000e+00}&2&9&0.01& $5$ & $30$  & $0.01$ & \num{ 2.000e+00} & \num{ 3.017e-06} \\
CRESC50 & $6$ & $100$ & $16$ & $22$  & $2.49$ & \num{ 6.057e-01} & \num{ 9.032e-11}&344&1534&2.38& $10$ & $91$  & $0.06$ & \num{ 6.567e-01} & \num{ 9.586e-04} \\
DEMBO7 & $16$ & $20$ & $13$ & $29$  & $1.44$ & \num{ 1.879e+02} & \num{ 5.917e-08}&198&512&1.19& $16$ & $358$  & $0.12$ & \num{ 1.748e+02} & \num{ 1.074e-05} \\
DNIEPER & $57$ & $24$ & $4$ & $61$  & $1.86$ & \num{ 1.941e+04} & \num{ 4.264e-11}&96&460&1.35& $5$ & $359$  & $0.07$ & \num{ 1.874e+04} & \num{ 3.740e-10} \\
EXPFITA & $5$ & $22$ & $1$ & $10999$  & $0.62$ & \num{ 2.999e+01}(***) & \num{ 1.798e+308}&0&2&0.0& $7$ & $56$  & $0.02$ & \num{ 9.889e-02} & \num{ 1.354e-14} \\
HIMMELBI & $100$ & $12$ & $114$ & $215$  & $70.51$ & \num{-1.467e+03}(*) & \num{ 1.789e-09}&2618&12588&68.05& $3$ & $408$  & $0.05$ & \num{-1.625e+03} & \num{ 0.000e+00} \\
SYNTHES1 & $6$ & $6$ & $1$ & $2999$  & $0.23$ & \num{ 1.000e+01}(***) & \num{ 1.798e+308}&0&2&0.0& $7$ & $65$  & $0.02$ & \num{ 7.593e-01} & \num{ 3.261e-07} \\
TWOBARS & $2$ & $2$ & $2$ & $4$  & $0.07$ & \num{ 1.523e+00} & \num{ 7.885e-09}&13&43&0.04& $4$ & $26$  & $0.01$ & \num{ 1.558e+00} & \num{ 0.000e+00} \\
DIPIGRI & $7$ & $4$ & $1$ & $8$  & $0.15$ & \num{ 7.140e+02} & \num{ 0.000e+00}&22&80&0.13& $0$ & $5$  & $0.0$ & \num{ 7.140e+02} & \num{ 0.000e+00} \\
\end{tabular}}
\caption{Noiseless Problems with Infeasible $x_0$. $\tau = 10^{-1}$; $n$: number of variables, $m$: number of constraints, \#iter: number of (outer) iterations, \#feval: number of function evaluations, time: total CPU time passed,  $f$: final objective value, feas err: final feasibility error, \#iter(sub): total number of iterations for solving TR subproblem, \#ceval: total number of constraint evaluations(note that the number of constraint evaluations and function evaluations are same for \texttt{KNITRO}), time(sub): total CPU time elapsed for solving subproblem. * indicates that \texttt{FIBO} terminates with singular interpolation system error, ** indicates that \texttt{FIBO} terminates with maximum number of iterations, *** indicates that \texttt{FIBO} terminates with maximum number of function evaluations. }
\label{tab:noiseless infeas 1e-1}
\end{table}

%tau = 1e-3
\begin{table}[htp]
\centering
\resizebox{0.77\columnwidth}{!}{\begin{tabular}{lll | llllllll | lllll}
\hline
\multicolumn{3}{l}{} & \multicolumn{8}{c}{FIBO}       & \multicolumn{5}{c}{knitro}         \\ \hline
Problem   & n   & m  & \#iter & \#feval & time & $f$ & feas err & \#iter(sub) & \#ceval & time(sub) &  \#iter & \#feval & time & $f$ & feas err \\ \hline
HS13 & $2$ & $1$ & $2$ & $4$  & $0.19$ & \num{ 1.000e+00} & \num{ 1.310e-14}&75&106&0.17& $18$ & $90$  & $0.08$ & \num{ 1.001e+00} & \num{ 0.000e+00} \\
HS22 & $2$ & $2$ & $1$ & $3$  & $0.04$ & \num{ 1.000e+00} & \num{ 1.789e-12}&9&25&0.02& $3$ & $16$  & $0.01$ & \num{ 1.000e+00} & \num{ 1.231e-06} \\
HS23 & $2$ & $5$ & $3$ & $5$  & $0.05$ & \num{ 2.000e+00} & \num{ 0.000e+00}&11&16&0.03& $4$ & $20$  & $0.01$ & \num{ 2.001e+00} & \num{ 0.000e+00} \\
HS26 & $3$ & $1$ & $36$ & $39$  & $1.86$ & \num{ 1.737e-04} & \num{ 1.040e-12}&754&4339&1.77& $14$ & $87$  & $0.03$ & \num{ 4.252e-07} & \num{ 9.633e-07} \\
HS32 & $3$ & $2$ & $4$ & $7$  & $0.06$ & \num{ 1.000e+00} & \num{ 0.000e+00}&10&17&0.03& $2$ & $15$  & $0.01$ & \num{ 1.000e+00} & \num{ 0.000e+00} \\
HS34 & $3$ & $2$ & $6$ & $9$  & $0.11$ & \num{-8.340e-01} & \num{ 1.776e-15}&32&71&0.08& $6$ & $39$  & $0.01$ & \num{-8.340e-01} & \num{ 5.446e-09} \\
HS40 & $4$ & $3$ & $2$ & $6$  & $0.05$ & \num{-2.491e-01} & \num{ 2.189e-09}&8&25&0.02& $4$ & $33$  & $0.01$ & \num{-2.500e-01} & \num{ 5.162e-08} \\
HS44 & $4$ & $6$ & $3$ & $8$  & $0.05$ & \num{-1.000e+00}(*) & \num{ 0.000e+00}&12&24&0.03& $4$ & $30$  & $0.01$ & \num{-1.500e+01} & \num{ 0.000e+00} \\
HS47 & $5$ & $3$ & $23$ & $28$  & $3.41$ & \num{ 1.989e-07} & \num{ 3.781e-11}&1391&7878&3.33& $72$ & $696$  & $0.36$ & \num{ 1.007e+01} & \num{ 6.661e-16} \\
HS50 & $5$ & $3$ & $19$ & $24$  & $1.02$ & \num{ 1.056e+00} & \num{ 2.398e-14}&126&423&0.96& $11$ & $85$  & $0.03$ & \num{ 1.339e-06} & \num{ 8.882e-16} \\
HS64 & $3$ & $1$ & $13$ & $16$  & $0.5$ & \num{ 6.303e+03} & \num{ 1.443e-15}&126&411&0.44& $9$ & $54$  & $0.02$ & \num{ 6.306e+03} & \num{ 0.000e+00} \\
HS66 & $3$ & $2$ & $2$ & $5$  & $0.03$ & \num{ 5.182e-01} & \num{ 3.563e-12}&7&13&0.02& $3$ & $21$  & $0.01$ & \num{ 5.186e-01} & \num{ 0.000e+00} \\
HS67 & $3$ & $14$ & $38$ & $41$  & $1.51$ & \num{-1.162e+03} & \num{ 0.000e+00}&448&1378&1.42& $9$ & $50$  & $0.02$ & \num{-1.161e+03} & \num{ 0.000e+00} \\
HS72 & $4$ & $2$ & $5$ & $9$  & $0.2$ & \num{ 7.284e+02} & \num{ 0.000e+00}&38&102&0.16& $11$ & $72$  & $0.05$ & \num{ 7.277e+02} & \num{ 8.766e-07} \\
HS75 & $4$ & $5$ & $16$ & $20$  & $0.35$ & \num{ 5.174e+03} & \num{ 5.684e-14}&97&185&0.28& $5$ & $36$  & $0.03$ & \num{ 5.174e+03} & \num{ 1.880e-04} \\
HS85 & $5$ & $21$ & $31$ & $36$  & $0.96$ & \num{-2.216e+00} & \num{ 0.000e+00}&357&904&0.87& $10$ & $78$  & $0.03$ & \num{-2.216e+00} & \num{ 2.389e-08} \\
HS87 & $6$ & $4$ & $1$ & $7$  & $0.09$ & \num{ 8.997e+03} & \num{ 2.274e-13}&9&27&0.06& $3$ & $33$  & $0.01$ & \num{ 8.997e+03} & \num{ 3.775e-08} \\
HS88 & $2$ & $1$ & $2$ & $4$  & $0.23$ & \num{ 1.363e+00} & \num{ 9.216e-16}&41&150&0.19& $16$ & $94$  & $0.05$ & \num{ 1.363e+00} & \num{ 1.074e-08} \\
HS89 & $3$ & $1$ & $4$ & $7$  & $0.6$ & \num{ 1.363e+00} & \num{ 3.663e-16}&101&518&0.48& $23$ & $147$  & $0.07$ & \num{ 1.362e+00} & \num{ 7.400e-07} \\
HS90 & $4$ & $1$ & $6$ & $10$  & $1.57$ & \num{ 1.363e+00} & \num{ 5.260e-17}&277&1551&1.52& $36$ & $334$  & $0.16$ & \num{ 1.363e+00} & \num{ 1.142e-07} \\
HS93 & $6$ & $2$ & $21$ & $27$  & $1.51$ & \num{ 1.352e+02} & \num{ 4.695e-11}&333&1707&1.43& $8$ & $87$  & $0.03$ & \num{ 1.351e+02} & \num{ 0.000e+00} \\
HS98 & $6$ & $4$ & $2$ & $8$  & $0.06$ & \num{ 3.136e+00} & \num{ 0.000e+00}&4&11&0.01& $6$ & $56$  & $0.01$ & \num{ 3.136e+00} & \num{ 0.000e+00} \\
HS100 & $7$ & $4$ & $51$ & $58$  & $15.19$ & \num{ 6.813e+02} & \num{ 1.388e-17}&1609&11183&15.01& $3$ & $55$  & $0.02$ & \num{ 6.811e+02} & \num{ 0.000e+00} \\
HS101 & $7$ & $5$ & $44$ & $51$  & $7.84$ & \num{ 1.812e+03} & \num{ 3.006e-10}&1239&7252&7.65& $26$ & $269$  & $0.1$ & \num{ 1.810e+03} & \num{ 2.249e-05} \\
HS102 & $7$ & $5$ & $69$ & $76$  & $10.15$ & \num{ 9.125e+02} & \num{ 4.272e-13}&1483&7671&9.89& $26$ & $282$  & $0.09$ & \num{ 9.120e+02} & \num{ 3.705e-06} \\
HS103 & $7$ & $5$ & $37$ & $44$  & $3.97$ & \num{ 5.437e+02} & \num{ 7.980e-16}&671&3472&3.8& $20$ & $212$  & $0.07$ & \num{ 5.434e+02} & \num{ 3.104e-04} \\
HS104 & $8$ & $5$ & $47$ & $55$  & $6.2$ & \num{ 3.952e+00} & \num{ 5.218e-15}&1210&6053&5.99& $11$ & $130$  & $0.03$ & \num{ 3.951e+00} & \num{ 5.769e-07} \\
LOADBAL & $31$ & $31$ & $1$ & $15499$  & $0.92$ & \num{ 1.547e+00}(***) & \num{ 1.798e+308}&0&2&0.0& $9$ & $330$  & $0.05$ & \num{ 4.531e-01} & \num{ 1.421e-14} \\
OPTPRLOC & $30$ & $30$ & $6$ & $36$  & $0.45$ & \num{-1.642e+01} & \num{ 0.000e+00}&47&114&0.39& $12$ & $456$  & $0.06$ & \num{-1.642e+01} & \num{ 1.279e-13} \\
CB3 & $3$ & $3$ & $1$ & $4$  & $0.03$ & \num{ 2.000e+00} & \num{ 0.000e+00}&2&9&0.01& $5$ & $30$  & $0.01$ & \num{ 2.000e+00} & \num{ 3.017e-06} \\
CRESC50 & $6$ & $100$ & $68$ & $74$  & $18.03$ & \num{ 5.942e-01} & \num{ 0.000e+00}&2416&11263&17.76& $456$ & $6023$  & $2.47$ & \num{ 5.948e-01} & \num{ 4.547e-13} \\
DEMBO7 & $16$ & $20$ & $18$ & $34$  & $2.17$ & \num{ 1.749e+02} & \num{ 1.067e-10}&301&798&2.03& $16$ & $358$  & $0.11$ & \num{ 1.748e+02} & \num{ 1.074e-05} \\
DNIEPER & $57$ & $24$ & $5$ & $62$  & $1.45$ & \num{ 1.874e+04} & \num{ 1.405e-09}&100&468&1.35& $5$ & $359$  & $0.07$ & \num{ 1.874e+04} & \num{ 3.740e-10} \\
EXPFITA & $5$ & $22$ & $1$ & $10999$  & $0.59$ & \num{ 2.999e+01}(***) & \num{ 1.798e+308}&0&2&0.0& $12$ & $91$  & $0.03$ & \num{ 1.355e-03} & \num{ 6.750e-14} \\
HIMMELBI & $100$ & $12$ & $114$ & $215$  & $68.61$ & \num{-1.467e+03}(*) & \num{ 1.789e-09}&2618&12588&67.1& $16$ & $1734$  & $0.19$ & \num{-1.734e+03} & \num{ 3.553e-15} \\
SYNTHES1 & $6$ & $6$ & $1$ & $2999$  & $0.22$ & \num{ 1.000e+01}(***) & \num{ 1.798e+308}&0&2&0.0& $7$ & $65$  & $0.02$ & \num{ 7.593e-01} & \num{ 3.261e-07} \\
TWOBARS & $2$ & $2$ & $3$ & $5$  & $0.08$ & \num{ 1.509e+00} & \num{ 4.742e-11}&16&60&0.05& $9$ & $49$  & $0.02$ & \num{ 1.509e+00} & \num{ 8.844e-08} \\
DIPIGRI & $7$ & $4$ & $51$ & $58$  & $11.34$ & \num{ 6.812e+02} & \num{ 6.939e-18}&1491&9510&11.16& $3$ & $55$  & $0.01$ & \num{ 6.811e+02} & \num{ 0.000e+00} \\
\end{tabular}}
\caption{Noiseless Problems with Infeasible $x_0$. $\tau = 10^{-3}$;  $n$: number of variables, $m$: number of constraints, \#iter: number of (outer) iterations, \#feval: number of function evaluations, time: total CPU time passed,  $f$: final objective value, feas err: final feasibility error, \#iter(sub): total number of iterations for solving TR subproblem, \#ceval: total number of constraint evaluations(note that the number of constraint evaluations and function evaluations are same for \texttt{KNITRO}), time(sub): total CPU time elapsed for solving subproblem. * indicates that \texttt{FIBO} terminates with singular interpolation system error, ** indicates that \texttt{FIBO} terminates with maximum number of iterations, *** indicates that \texttt{FIBO} terminates with maximum number of function evaluations. }
\label{tab:noiseless infeas 1e-3}
\end{table}

%tau = 1e-5
\begin{table}[htp]
\centering
\resizebox{0.77\columnwidth}{!}{\begin{tabular}{lll | llllllll | lllll}
\hline
\multicolumn{3}{l}{} & \multicolumn{8}{c}{FIBO}       & \multicolumn{5}{c}{knitro}         \\ \hline
Problem   & n   & m  & \#iter & \#feval & time & $f$ & feas err & \#iter(sub) & \#ceval & time(sub) &  \#iter & \#feval & time & $f$ & feas err \\ \hline
HS13 & $2$ & $1$ & $3$ & $5$  & $0.29$ & \num{ 1.000e+00} & \num{ 1.309e-14}&136&191&0.28& $19$ & $102$  & $0.09$ & \num{ 9.983e-01} & \num{ 6.579e-10} \\
HS22 & $2$ & $2$ & $1$ & $3$  & $0.04$ & \num{ 1.000e+00} & \num{ 1.789e-12}&9&25&0.02& $3$ & $16$  & $0.01$ & \num{ 1.000e+00} & \num{ 1.231e-06} \\
HS23 & $2$ & $5$ & $3$ & $5$  & $0.04$ & \num{ 2.000e+00} & \num{ 0.000e+00}&11&16&0.02& $5$ & $24$  & $0.01$ & \num{ 2.000e+00} & \num{ 0.000e+00} \\
HS26 & $3$ & $1$ & $38$ & $41$  & $8.15$ & \num{ 8.872e-07} & \num{ 4.441e-16}&3356&26977&8.07& $14$ & $87$  & $0.03$ & \num{ 4.252e-07} & \num{ 9.633e-07} \\
HS32 & $3$ & $2$ & $4$ & $7$  & $0.05$ & \num{ 1.000e+00} & \num{ 0.000e+00}&10&17&0.03& $2$ & $15$  & $0.01$ & \num{ 1.000e+00} & \num{ 0.000e+00} \\
HS34 & $3$ & $2$ & $6$ & $9$  & $0.11$ & \num{-8.340e-01} & \num{ 1.776e-15}&32&71&0.08& $6$ & $39$  & $0.01$ & \num{-8.340e-01} & \num{ 5.446e-09} \\
HS40 & $4$ & $3$ & $4$ & $8$  & $0.08$ & \num{-2.500e-01} & \num{ 1.058e-12}&13&36&0.04& $4$ & $33$  & $0.01$ & \num{-2.500e-01} & \num{ 5.162e-08} \\
HS44 & $4$ & $6$ & $3$ & $8$  & $0.05$ & \num{-1.000e+00}(*) & \num{ 0.000e+00}&12&24&0.04& $4$ & $30$  & $0.01$ & \num{-1.500e+01} & \num{ 0.000e+00} \\
HS47 & $5$ & $3$ & $23$ & $28$  & $3.42$ & \num{ 1.989e-07} & \num{ 3.781e-11}&1391&7878&3.35& $72$ & $696$  & $0.36$ & \num{ 1.007e+01} & \num{ 6.661e-16} \\
HS50 & $5$ & $3$ & $19$ & $24$  & $1.11$ & \num{ 1.056e+00} & \num{ 2.398e-14}&126&423&1.06& $11$ & $85$  & $0.03$ & \num{ 1.339e-06} & \num{ 8.882e-16} \\
HS64 & $3$ & $1$ & $17$ & $20$  & $0.58$ & \num{ 6.300e+03} & \num{ 0.000e+00}&156&447&0.52& $11$ & $64$  & $0.03$ & \num{ 6.300e+03} & \num{ 6.916e-06} \\
HS66 & $3$ & $2$ & $2$ & $5$  & $0.03$ & \num{ 5.182e-01} & \num{ 3.563e-12}&7&13&0.02& $4$ & $27$  & $0.01$ & \num{ 5.182e-01} & \num{ 0.000e+00} \\
HS67 & $3$ & $14$ & $39$ & $42$  & $1.46$ & \num{-1.162e+03} & \num{ 0.000e+00}&450&1381&1.38& $10$ & $55$  & $0.03$ & \num{-1.162e+03} & \num{ 0.000e+00} \\
HS72 & $4$ & $2$ & $6$ & $10$  & $0.27$ & \num{ 7.277e+02} & \num{ 6.263e-14}&43&109&0.22& $11$ & $72$  & $0.04$ & \num{ 7.277e+02} & \num{ 8.766e-07} \\
HS75 & $4$ & $5$ & $16$ & $20$  & $0.34$ & \num{ 5.174e+03} & \num{ 5.684e-14}&97&185&0.27& $5$ & $36$  & $0.03$ & \num{ 5.174e+03} & \num{ 1.880e-04} \\
HS85 & $5$ & $21$ & $31$ & $36$  & $0.94$ & \num{-2.216e+00} & \num{ 0.000e+00}&357&904&0.83& $10$ & $78$  & $0.03$ & \num{-2.216e+00} & \num{ 2.389e-08} \\
HS87 & $6$ & $4$ & $2$ & $8$  & $0.13$ & \num{ 8.997e+03} & \num{ 7.844e-12}&15&46&0.1& $3$ & $33$  & $0.02$ & \num{ 8.997e+03} & \num{ 3.775e-08} \\
HS88 & $2$ & $1$ & $4$ & $6$  & $0.34$ & \num{ 1.363e+00} & \num{ 1.301e-14}&69&229&0.3& $16$ & $94$  & $0.05$ & \num{ 1.363e+00} & \num{ 1.074e-08} \\
HS89 & $3$ & $1$ & $5$ & $8$  & $0.64$ & \num{ 1.363e+00} & \num{ 0.000e+00}&113&576&0.54& $23$ & $147$  & $0.08$ & \num{ 1.362e+00} & \num{ 7.400e-07} \\
HS90 & $4$ & $1$ & $11$ & $15$  & $3.63$ & \num{ 1.363e+00} & \num{ 8.118e-18}&718&3460&3.57& $36$ & $334$  & $0.18$ & \num{ 1.363e+00} & \num{ 1.142e-07} \\
HS93 & $6$ & $2$ & $25$ & $31$  & $1.58$ & \num{ 1.351e+02} & \num{ 3.331e-16}&371&1843&1.5& $8$ & $87$  & $0.03$ & \num{ 1.351e+02} & \num{ 0.000e+00} \\
HS98 & $6$ & $4$ & $2$ & $8$  & $0.05$ & \num{ 3.136e+00} & \num{ 0.000e+00}&4&11&0.02& $6$ & $56$  & $0.01$ & \num{ 3.136e+00} & \num{ 0.000e+00} \\
HS100 & $7$ & $4$ & $89$ & $96$  & $24.87$ & \num{ 6.806e+02} & \num{ 0.000e+00}&2860&19400&24.54& $14$ & $212$  & $0.07$ & \num{ 6.806e+02} & \num{ 0.000e+00} \\
HS101 & $7$ & $5$ & $71$ & $78$  & $12.77$ & \num{ 1.810e+03} & \num{ 2.776e-16}&2004&11927&12.48& $27$ & $280$  & $0.1$ & \num{ 1.810e+03} & \num{ 2.398e-05} \\
HS102 & $7$ & $5$ & $78$ & $85$  & $10.48$ & \num{ 9.119e+02} & \num{ 0.000e+00}&1601&8233&10.11& $28$ & $305$  & $0.09$ & \num{ 9.119e+02} & \num{ 2.176e-06} \\
HS103 & $7$ & $5$ & $39$ & $46$  & $4.26$ & \num{ 5.437e+02} & \num{ 6.106e-16}&685&3522&4.09& $20$ & $212$  & $0.07$ & \num{ 5.434e+02} & \num{ 3.104e-04} \\
HS104 & $8$ & $5$ & $80$ & $88$  & $13.33$ & \num{ 3.951e+00} & \num{ 0.000e+00}&2262&12387&13.0& $11$ & $130$  & $0.03$ & \num{ 3.951e+00} & \num{ 5.769e-07} \\
LOADBAL & $31$ & $31$ & $1$ & $15499$  & $0.89$ & \num{ 1.547e+00}(***) & \num{ 1.798e+308}&0&2&0.0& $17$ & $594$  & $0.09$ & \num{ 4.529e-01} & \num{ 1.483e-14} \\
OPTPRLOC & $30$ & $30$ & $8$ & $38$  & $0.54$ & \num{-1.642e+01} & \num{ 4.504e-10}&53&128&0.44& $12$ & $456$  & $0.06$ & \num{-1.642e+01} & \num{ 1.279e-13} \\
CB3 & $3$ & $3$ & $1$ & $4$  & $0.03$ & \num{ 2.000e+00} & \num{ 0.000e+00}&2&9&0.01& $5$ & $30$  & $0.01$ & \num{ 2.000e+00} & \num{ 3.017e-06} \\
CRESC50 & $6$ & $100$ & $1000$ & $1006$  & $214.67$ & \num{ 5.939e-01}(**) & \num{ 5.482e-13}&26160&96697&211.17& $456$ & $6023$  & $2.37$ & \num{ 5.948e-01} & \num{ 4.547e-13} \\
DEMBO7 & $16$ & $20$ & $19$ & $35$  & $2.11$ & \num{ 1.748e+02} & \num{ 4.073e-15}&302&800&1.97& $16$ & $358$  & $0.1$ & \num{ 1.748e+02} & \num{ 1.074e-05} \\
DNIEPER & $57$ & $24$ & $5$ & $62$  & $1.43$ & \num{ 1.874e+04} & \num{ 1.405e-09}&100&468&1.36& $5$ & $359$  & $0.07$ & \num{ 1.874e+04} & \num{ 3.740e-10} \\
EXPFITA & $5$ & $22$ & $1$ & $10999$  & $0.63$ & \num{ 2.999e+01}(***) & \num{ 1.798e+308}&0&2&0.0& $13$ & $98$  & $0.03$ & \num{ 1.137e-03} & \num{ 6.047e-12} \\
HIMMELBI & $100$ & $12$ & $114$ & $215$  & $68.54$ & \num{-1.467e+03}(*) & \num{ 1.789e-09}&2618&12588&67.29& $41$ & $4291$  & $0.45$ & \num{-1.736e+03} & \num{ 1.421e-14} \\
SYNTHES1 & $6$ & $6$ & $1$ & $2999$  & $0.21$ & \num{ 1.000e+01}(***) & \num{ 1.798e+308}&0&2&0.0& $7$ & $65$  & $0.01$ & \num{ 7.593e-01} & \num{ 3.261e-07} \\
TWOBARS & $2$ & $2$ & $3$ & $5$  & $0.08$ & \num{ 1.509e+00} & \num{ 4.742e-11}&16&60&0.05& $9$ & $49$  & $0.02$ & \num{ 1.509e+00} & \num{ 8.844e-08} \\
DIPIGRI & $7$ & $4$ & $87$ & $94$  & $78.71$ & \num{ 6.806e+02} & \num{ 1.265e-13}&7019&23031&78.4& $14$ & $212$  & $0.06$ & \num{ 6.806e+02} & \num{ 0.000e+00} \\
\end{tabular}}
\caption{Noiseless Problems with Infeasible $x_0$ $\tau = 10^{-5}$;  $n$: number of variables, $m$: number of constraints, \#iter: number of (outer) iterations, \#feval: number of function evaluations, time: total CPU time passed,  $f$: final objective value, feas err: final feasibility error, \#iter(sub): total number of iterations for solving TR subproblem, \#ceval: total number of constraint evaluations(note that the number of constraint evaluations and function evaluations are same for \texttt{KNITRO}), time(sub): total CPU time elapsed for solving subproblem. * indicates that \texttt{FIBO} terminates with singular interpolation system error, ** indicates that \texttt{FIBO} terminates with maximum number of iterations, *** indicates that \texttt{FIBO} terminates with maximum number of function evaluations. }
\label{tab:noiseless infeas 1e-5}
\end{table}

%tau = 1e-7
\begin{table}[htp]
\centering
\resizebox{0.77\columnwidth}{!}{\begin{tabular}{lll | llllllll | lllll}
\hline
\multicolumn{3}{l}{} & \multicolumn{8}{c}{FIBO}       & \multicolumn{5}{c}{knitro}         \\ \hline
Problem   & n   & m  & \#iter & \#feval & time & $f$ & feas err & \#iter(sub) & \#ceval & time(sub) &  \#iter & \#feval & time & $f$ & feas err \\ \hline
HS13 & $2$ & $1$ & $3$ & $5$  & $0.29$ & \num{ 1.000e+00} & \num{ 1.309e-14}&136&191&0.28& $19$ & $102$  & $0.09$ & \num{ 9.983e-01} & \num{ 6.579e-10} \\
HS22 & $2$ & $2$ & $1$ & $3$  & $0.04$ & \num{ 1.000e+00} & \num{ 1.789e-12}&9&25&0.02& $3$ & $16$  & $0.01$ & \num{ 1.000e+00} & \num{ 1.231e-06} \\
HS23 & $2$ & $5$ & $3$ & $5$  & $0.05$ & \num{ 2.000e+00} & \num{ 0.000e+00}&11&16&0.03& $6$ & $28$  & $0.01$ & \num{ 2.000e+00} & \num{ 0.000e+00} \\
HS26 & $3$ & $1$ & $43$ & $46$  & $12.41$ & \num{ 2.878e-08} & \num{ 1.110e-16}&5194&38851&12.32& $28$ & $203$  & $0.08$ & \num{ 2.719e-11} & \num{ 8.518e-07} \\
HS32 & $3$ & $2$ & $4$ & $7$  & $0.05$ & \num{ 1.000e+00} & \num{ 0.000e+00}&10&17&0.03& $2$ & $15$  & $0.01$ & \num{ 1.000e+00} & \num{ 0.000e+00} \\
HS34 & $3$ & $2$ & $6$ & $9$  & $0.09$ & \num{-8.340e-01} & \num{ 1.776e-15}&32&71&0.07& $6$ & $39$  & $0.01$ & \num{-8.340e-01} & \num{ 5.446e-09} \\
HS40 & $4$ & $3$ & $4$ & $8$  & $0.07$ & \num{-2.500e-01} & \num{ 1.058e-12}&13&36&0.04& $4$ & $33$  & $0.01$ & \num{-2.500e-01} & \num{ 5.162e-08} \\
HS44 & $4$ & $6$ & $3$ & $8$  & $0.05$ & \num{-1.000e+00}(*) & \num{ 0.000e+00}&12&24&0.03& $4$ & $30$  & $0.01$ & \num{-1.500e+01} & \num{ 0.000e+00} \\
HS47 & $5$ & $3$ & $23$ & $28$  & $3.44$ & \num{ 1.989e-07} & \num{ 3.781e-11}&1391&7878&3.36& $72$ & $696$  & $0.33$ & \num{ 1.007e+01} & \num{ 6.661e-16} \\
HS50 & $5$ & $3$ & $19$ & $24$  & $1.12$ & \num{ 1.056e+00} & \num{ 2.398e-14}&126&423&1.06& $13$ & $100$  & $0.03$ & \num{ 2.358e-08} & \num{ 4.441e-16} \\
HS64 & $3$ & $1$ & $17$ & $20$  & $0.59$ & \num{ 6.300e+03} & \num{ 0.000e+00}&156&447&0.52& $11$ & $64$  & $0.03$ & \num{ 6.300e+03} & \num{ 6.916e-06} \\
HS66 & $3$ & $2$ & $2$ & $5$  & $0.03$ & \num{ 5.182e-01} & \num{ 3.563e-12}&7&13&0.02& $5$ & $33$  & $0.01$ & \num{ 5.182e-01} & \num{ 5.154e-11} \\
HS67 & $3$ & $14$ & $40$ & $43$  & $1.46$ & \num{-1.162e+03} & \num{ 2.274e-12}&458&1390&1.36& $11$ & $60$  & $0.03$ & \num{-1.162e+03} & \num{ 1.364e-12} \\
HS72 & $4$ & $2$ & $7$ & $11$  & $0.25$ & \num{ 7.277e+02} & \num{ 2.819e-18}&44&111&0.2& $11$ & $72$  & $0.04$ & \num{ 7.277e+02} & \num{ 8.766e-07} \\
HS75 & $4$ & $5$ & $16$ & $20$  & $0.35$ & \num{ 5.174e+03} & \num{ 5.684e-14}&97&185&0.28& $5$ & $36$  & $0.03$ & \num{ 5.174e+03} & \num{ 1.880e-04} \\
HS85 & $5$ & $21$ & $31$ & $36$  & $1.42$ & \num{-2.216e+00} & \num{ 0.000e+00}&357&904&1.3& $10$ & $78$  & $0.03$ & \num{-2.216e+00} & \num{ 2.389e-08} \\
HS87 & $6$ & $4$ & $9$ & $15$  & $0.7$ & \num{ 8.997e+03} & \num{ 2.274e-13}&80&287&0.67& $6$ & $75$  & $0.02$ & \num{ 8.997e+03} & \num{ 7.329e-11} \\
HS88 & $2$ & $1$ & $4$ & $6$  & $0.34$ & \num{ 1.363e+00} & \num{ 1.301e-14}&69&229&0.3& $16$ & $94$  & $0.05$ & \num{ 1.363e+00} & \num{ 1.074e-08} \\
HS89 & $3$ & $1$ & $6$ & $9$  & $0.78$ & \num{ 1.363e+00} & \num{ 5.535e-16}&142&676&0.66& $23$ & $147$  & $0.07$ & \num{ 1.362e+00} & \num{ 7.400e-07} \\
HS90 & $4$ & $1$ & $11$ & $15$  & $3.63$ & \num{ 1.363e+00} & \num{ 8.118e-18}&718&3460&3.57& $36$ & $334$  & $0.18$ & \num{ 1.363e+00} & \num{ 1.142e-07} \\
HS93 & $6$ & $2$ & $36$ & $42$  & $2.19$ & \num{ 1.351e+02} & \num{ 0.000e+00}&483&2292&2.06& $9$ & $96$  & $0.03$ & \num{ 1.351e+02} & \num{ 0.000e+00} \\
HS98 & $6$ & $4$ & $2$ & $8$  & $0.05$ & \num{ 3.136e+00} & \num{ 0.000e+00}&4&11&0.02& $6$ & $56$  & $0.01$ & \num{ 3.136e+00} & \num{ 0.000e+00} \\
HS100 & $7$ & $4$ & $106$ & $113$  & $28.46$ & \num{ 6.806e+02} & \num{ 0.000e+00}&3280&22448&28.08& $28$ & $362$  & $0.11$ & \num{ 6.806e+02} & \num{ 5.093e-11} \\
HS101 & $7$ & $5$ & $90$ & $97$  & $15.51$ & \num{ 1.810e+03} & \num{ 2.742e-12}&2377&14126&15.15& $28$ & $289$  & $0.11$ & \num{ 1.810e+03} & \num{ 1.302e-05} \\
HS102 & $7$ & $5$ & $89$ & $96$  & $12.12$ & \num{ 9.119e+02} & \num{ 2.012e-16}&1847&9687&11.74& $29$ & $314$  & $0.1$ & \num{ 9.119e+02} & \num{ 7.385e-07} \\
HS103 & $7$ & $5$ & $42$ & $49$  & $4.11$ & \num{ 5.437e+02} & \num{ 4.302e-16}&713&3625&3.92& $20$ & $212$  & $0.07$ & \num{ 5.434e+02} & \num{ 3.104e-04} \\
HS104 & $8$ & $5$ & $88$ & $96$  & $13.73$ & \num{ 3.951e+00} & \num{ 4.441e-16}&2400&13131&13.38& $11$ & $130$  & $0.03$ & \num{ 3.951e+00} & \num{ 5.769e-07} \\
LOADBAL & $31$ & $31$ & $1$ & $15499$  & $0.9$ & \num{ 1.547e+00}(***) & \num{ 1.798e+308}&0&2&0.0& $18$ & $627$  & $0.09$ & \num{ 4.529e-01} & \num{ 8.549e-15} \\
OPTPRLOC & $30$ & $30$ & $8$ & $38$  & $0.53$ & \num{-1.642e+01} & \num{ 4.504e-10}&53&128&0.43& $12$ & $456$  & $0.06$ & \num{-1.642e+01} & \num{ 1.279e-13} \\
CB3 & $3$ & $3$ & $1$ & $4$  & $0.02$ & \num{ 2.000e+00} & \num{ 0.000e+00}&2&9&0.0& $5$ & $30$  & $0.01$ & \num{ 2.000e+00} & \num{ 3.017e-06} \\
CRESC50 & $6$ & $100$ & $1000$ & $1006$  & $214.2$ & \num{ 5.939e-01}(**) & \num{ 5.482e-13}&26160&96697&210.8& $456$ & $6023$  & $2.34$ & \num{ 5.948e-01} & \num{ 4.547e-13} \\
DEMBO7 & $16$ & $20$ & $19$ & $35$  & $2.11$ & \num{ 1.748e+02} & \num{ 4.073e-15}&302&795&2.04& $16$ & $358$  & $0.11$ & \num{ 1.748e+02} & \num{ 1.074e-05} \\
DNIEPER & $57$ & $24$ & $23$ & $80$  & $3.47$ & \num{ 1.874e+04} & \num{ 5.684e-14}&272&1091&3.27& $5$ & $359$  & $0.06$ & \num{ 1.874e+04} & \num{ 3.740e-10} \\
EXPFITA & $5$ & $22$ & $1$ & $10999$  & $0.6$ & \num{ 2.999e+01}(***) & \num{ 1.798e+308}&0&2&0.0& $14$ & $105$  & $0.04$ & \num{ 1.137e-03} & \num{ 1.776e-14} \\
HIMMELBI & $100$ & $12$ & $114$ & $215$  & $69.13$ & \num{-1.467e+03}(*) & \num{ 1.789e-09}&2618&12588&68.01& $46$ & $4805$  & $0.48$ & \num{-1.736e+03} & \num{ 2.132e-14} \\
SYNTHES1 & $6$ & $6$ & $1$ & $2999$  & $0.24$ & \num{ 1.000e+01}(***) & \num{ 1.798e+308}&0&2&0.0& $7$ & $65$  & $0.02$ & \num{ 7.593e-01} & \num{ 3.261e-07} \\
TWOBARS & $2$ & $2$ & $5$ & $7$  & $0.11$ & \num{ 1.509e+00} & \num{ 0.000e+00}&21&89&0.08& $9$ & $49$  & $0.02$ & \num{ 1.509e+00} & \num{ 8.844e-08} \\
DIPIGRI & $7$ & $4$ & $96$ & $103$  & $80.89$ & \num{ 6.806e+02} & \num{ 4.547e-13}&7235&24606&80.56& $28$ & $365$  & $0.11$ & \num{ 6.806e+02} & \num{ 7.237e-11} \\
\end{tabular}}
\caption{Noiseless Problems with Infeasible $x_0$ $\tau = 10^{-7}$; $n$: number of variables, $m$: number of constraints, \#iter: number of (outer) iterations, \#feval: number of function evaluations, time: total CPU time passed,  $f$: final objective value, feas err: final feasibility error, \#iter(sub): total number of iterations for solving TR subproblem, \#ceval: total number of constraint evaluations(note that the number of constraint evaluations and function evaluations are same for \texttt{KNITRO}), time(sub): total CPU time elapsed for solving subproblem. * indicates that \texttt{FIBO} terminates with singular interpolation system error, ** indicates that \texttt{FIBO} terminates with maximum number of iterations, *** indicates that \texttt{FIBO} terminates with maximum number of function evaluations. }
\label{tab:noiseless infeas 1e-7}
\end{table}

\end{landscape}

\end{document}